# Risk-Based Admissibility Assessment of Wind Generation Integrated into a Bulk Power System

Cheng Wang, Feng Liu, *Member, IEEE*, Jianhui Wang, *Senior Member, IEEE*, Wei Wei, Shengwei Mei, *Fellow, IEEE*

***Abstract*—The increasing integration of large-scale volatile and uncertain wind generation has brought great challenges to power system operations. In this paper, a risk-based admissibility assessment approach is proposed to quantitatively evaluate how much wind generation can be accommodated by the bulk power system under a given unit commitment (UC) strategy. Firstly, the operational risk brought by the variation and uncertainty of wind generation is developed as an admissibility measure of wind generation. Then its linear approximation is derived for practical implementation. Furthermore, a risk-minimization model is established to mathematically characterize the admissible region of wind generation. This model can be solved effectively by a modified column and constraint generation (C&CG) algorithm. Simulations on the IEEE 9-bus system and the real Guangdong power grid demonstrate the effectiveness and efficiency of the proposed methodology.**

## NOMENCLATURE

*Indices*

| | |
|---|---|
| $g$ | Index for generators. |
| $m$ | Index for wind farms. |
| $l$ | Index for transmission lines. |
| $j$ | Index for loads. |
| $n$ | Index for nodes. |
| $t$ | Index for time periods. |

*Parameters*

| | |
|---|---|
| $T$ | Number of time periods. |
| $N$ | Number of nodes. |
| $M$ | Number of wind farms. |
| $G$ | Number of thermal generators. |
| $P_g^{min}$ | Minimal output of generator $g$ when it is on. |
| $P_g^{max}$ | Maximal output of generator $g$ when it is on. |
| $R_+^g/R_-^g$ | Ramp-up/down limit for generator $g$. |
| $u_{gt}$ | Binary variable indicating whether generator $g$ is on or off in period $t$. |
| $r$ | Spinning reserve rate. |
| $F_l$ | Transmission capacity of line $l$. |
| $W$ | Wind generation uncertainty set. |
| $\hat{w}_{mt}$ | Forecasted output of wind farm $m$ in period $t$. |
| $w_m^{max}$ | Installed capacity of wind farm $m$. |
| $\Gamma^S/\Gamma^T$ | Budget of uncertainties over spatial/time scale. |
| $D_{jt}$ | Load demand of load node $j$ in period $t$. |
| $B$ | Node admittance matrix of the grid. |
| $Line_l$ | Indices of initial node and terminal node of line $l$. |
| $o_1/o_2$ | Number of initial/ terminal node of line $l$. |
| $\Phi(n)$ | The set of nodes connecting to node $n$. |
| $\theta_{nt}$ | Phase angle of node $n$ in period $t$. |
| $\alpha_{mt}$ | Confidence level of wind generation output interval. |
| $\beta_t/\beta_s$ | Confidence level of $\Gamma^T$/ $\Gamma^S$. |
| $e_{mt}$ | Price of wind generation curtailment of wind farm $m$ in period $t$ |
| $f_{jt}$ | Price of load shedding of load bus $j$ in period $t$ |

*Decision Variables*

| | |
|---|---|
| $P_{gt}$ | Output of generator $g$ in period $t$. |
| $v_{mt}^u$/ $v_{mt}^l$ | Binary variable indicating normalized positive/negative output deviation of wind farm $m$ in period $t$. |
| $\Delta w_{mt}$ | Wind generation curtailment in wind farm $m$ in period $t$. |
| $\Delta D_{jt}$ | Load shedding at load node $j$ in period $t$. |
| $w_{mt}^u$ | Upper output bound of wind farm $m$ in period $t$. |
| $w_{mt}^l$ | Lower output bound of wind farm $m$ in period $t$. |
| $Q_{mt}^p$ $/Q_{mt}^n$ | Auxiliary variables representing operational risk due to underestimation/ overestimation of the output of wind farm $m$ in period $t$. |

## I. INTRODUCTION

The increasing integration of large-scale volatile and uncertain wind power generation into bulk power systems have created great challenges to power system operations, particularly unit commitment (UC) and economic dispatch (ED). Much work has been done to enhance the capability of power systems to admit high-penetration wind generation. The existing literature on wind power integration can be roughly divided into two categories: high-accuracy wind generation forecasting as well as uncertainty modeling on the wind farm side, and flexible dispatch strategies on the power system side.

According to the forecast time horizon [1], wind generation forecasts can be categorized to the immediate short-term forecast, the short-term forecast and the mid-term to long-term forecast, respectively. Generally, forecast accuracy deteriorates as the forecast time window increases [2]. Practical experience indicates that the day-ahead wind generation forecast error usually varies from 10% to 20%

This work was supported in part by the Special Grant from EPRI of China (XTB51201303968) and Foundation for Innovative Research Groups of the National Natural Science Foundation of China (51321005).

C. Wang, F. Liu, W. Wei, and S. Mei are with the State Key Laboratory of Power Systems, Department of Electrical Engineering and Applied Electronic Technology, Tsinghua University, 100084 Beijing, China. (Corresponding author: Feng Liu, e-mail: lfeng@mail.tsinghua.edu.cn).

J. Wang is with the Argonne National Laboratory, Argonne, IL 60439, USA (e-mail: jianhui.wang@anl.gov)

using current forecast techniques, which is much larger than that of conventional load forecasts, and the effect of forecast errors must be taken into account when making dispatch decisions. Thus, it is crucial to model the uncertainty of wind generation appropriately. Two major approaches have been proposed in the literature. One is the scenario-based approach, where the discrete distribution of a limited number of scenarios is generated to approximate the probability density function (PDF) of wind generation [3]. This approach has been extensively applied to the decision making of power system dispatch under uncertainties [4]. The other is the uncertainty set approach, where a set of inequalities are used to characterize the variation and uncertainty of wind generation, including the upper and lower bounds of uncertainty in each period, the intensity of power fluctuation, and the correlation of forecast errors of different wind farm outputs [5], [6]. It should be noted that, both kinds of wind generation uncertainty models are constructed based on a probability associated with certain confidence levels of interest explicitly or implicitly.

As for the dispatch decision-making, both real-time balancing dispatch and day-ahead scheduling are involved. Current research pays more attention to the latter, especially the UC decision-making problem, as it is essential for the operational flexibility of a power system in a day. In a traditional security-constrained unit commitment (SCUC) problem, a certain level of spinning reserve rate (SRR) is usually required to hedge against the risk of unexpected load variations or generator outages [7]. When large-scale wind generation integration is considered, a common practical solution to cope with wind power uncertainty is to directly scale up SRR in the traditional SCUC model. The latest research, however, evidently shows that this approach may not be accurate [8], [9]. To solve this problem, a stochastic unit commitment (SUC) formulation has been proposed in [10], [11], [12], [13]. However, the performance of SUC might not be guaranteed as certain rare scenarios with a high impact may not be included in the process of scenario generation. To circumvent such a problem, robust unit commitment (RUC) is developed as it can generate robust optimal strategies that ensure the feasibility of the solution for all possible scenarios within a given uncertainty set [14], [15]. With either the SUC or the RUC, the operational flexibility of a power system will be effectively enhanced, resulting in higher admissibility of a power system in integrating wind generation.

As mentioned above, the majority of existing models of wind generation uncertainty are established based on a probability associated with certain confidence levels, which may not cover all possible realizations of wind generation. As a result, power system operational feasibility may not be guaranteed even if the SUC or the RUC is adopted, as the uncertainty models may fail to cover certain rare albeit high-risk scenarios (such as a very high wind power ramp due to a sudden gust). When such scenarios happen, emergency measures such as load shedding, wind generation curtailment or additional fast reserve provision may have to be exercised to ensure the reliability of the system. This matter of fact gives rise to an important issue: how much volatile and uncertain wind generation can be admitted by the bulk power system at most, provided that a UC strategy and a wind generation forecast are given? In other words, given that the current system operation methods cannot handle all the wind power variation scenarios, an admissibility assessment of wind power penetration in a specific power system should be conducted to define the actual level of wind power whose fluctuations the system can accommodate. To answer this question, [16] proposes an inspiring concept of “Do-Not-Exceed” (DNE) wind power integration limit under a given day-ahead ED strategy as well as the associated computation method. It essentially divides the wind generation into two parts: the admissible region and the inadmissible region. Within the admissible region, any realization of wind generation can be fully admitted with no need of additional emergency resources, which is regarded as riskless. In the inadmissible region, however, additional regulation resources may have to be used to eliminate the power imbalance due to unforeseen variations of wind generation. In this paper, the expected cost for such additional emergency regulations to recover the feasibility is referred to as *operational risk*. Inspired by the work in [16], this paper proposes an operational risk based assessment methodology to quantify the maximum admissible wind generation under a given UC strategy. Compared with existing work, the main contributions of this paper are threefold.

1) In [16], conventional generators are divided into two categories: the corrective control units (CCUs) and the non-CCUs. DNE limit is derived under a given ED strategy for both the CCUs and the non-CCUs. The pre-given ED decision, however, may not be optimal for the accommodation of wind generation. Consequently, this may result in underestimation of the dispatch capability of the power system to admit wind generation. In this regard, this paper alternatively studies the admissible region of wind generation under a given UC strategy without any restriction on ED strategies, which differs from the work in [16].

2) As the admissibility assessment is essentially a dynamic programming problem involving a multi-period decision-making process, the ED strategy in a certain period will exert impacts on the admissibility of wind generation afterwards. Thus, different ED strategies will result in different admissible wind generation. In this paper, the operational risk is taken as an admissibility measure so that the admissible region can be determined in the sense of minimum operational risk created by the inadmissible wind generation. Compared with the work [16], this approach takes the information of wind generation forecast errors into account, and provides more insights on the admissibility of wind generation.

3) Mathematically, the operational risk based admissibility assessment problem leads to a two-stage optimization model that can be directly solved by robust optimization algorithms [8], [14], [15]. To further improve the computing efficiency, an equivalent model is derived by imposing a penalty function on the original objective function of the master problem. Based on this, the column-and-constraint generation (C&CG) algorithm presented in [17] can be modified to solve the admissibility problem efficiently.

The remaining part of the paper is organized as follows. Section II describes the mathematical formulation. Section III presents the solution methodology. Section IV gives an illustrative example for the proposed model and algorithm. Section V demonstrates the intended applications of our work by using the data of a real provincial power grid in China. Finally, section VI concludes the paper with some discussion. In the Appendix, Some details about the mathematical formulation of approximation of operational risk is presented.

## II. Mathematical Formulation

### A. Admissibility of Wind Generation

Theoretically, when a UC strategy is given, the maximum admissible wind generation is essentially determined. Here the "admissible wind generation" is referred to as a subset of the wind generation uncertainty set, within which any scenario of wind generation will never cause operational infeasibility. As for admissible wind generation, no load shedding or wind generation curtailment is needed for reliability reasons. In this regard, the wind generation admissibility can be examined by solving the following relaxed bi-level program.

$$\max_{v^u, v^l} \min_{p,\theta,\Delta w,\Delta D} F = \sum_{t=1}^{T}\left(\sum_{m=1}^{M} e_{mt}\Delta w_{mt} + \sum_{j=1}^{J} f_{jt}\Delta D_{jt}\right) \quad (1a)$$

$$s.t. \quad u_{gt}P_{\min}^{g} \le p_{gt} \le u_{gt}P_{\max}^{g} \ \ \forall g, \forall t \quad (1b)$$

$$p_{gt} - p_{g(t+1)} \le u_{g(t+1)}R_{-}^{g} + (1-u_{g(t+1)})P_{\max}^{g} \quad \forall g, \ \forall t \quad (1c)$$

$$p_{g(t+1)} - p_{gt} \le u_{gt}R_{+}^{g} + (1-u_{gt})P_{\max}^{g} \quad \forall g, \ \forall t \quad (1d)$$

$$-F_l \le B_{o_1o_2}(\theta_{o_1t} - \theta_{o_2t}) \le F_l \quad o_1, o_2 \in Line_l, \forall l, \ \forall t \quad (1e)$$

$$-\pi \le \theta_{nt} \le \pi \quad \forall n, \forall t \quad (1f)$$

$$\theta_{reft} = 0 \quad \forall t \quad (1g)$$

$$\sum_{g\in\phi(n)} p_{gt} + \sum_{m\in\phi(n)} (w_{mt} - \Delta w_{mt}) - \sum_{o\in\phi(n)} B_{on}(\theta_{ot} - \theta_{ot}) - \sum_{j\in\phi(n)} (D_{jt} - \Delta D_{jt}) = 0 \quad \forall n, \ \forall t \quad (1h)$$

$$0 \le \Delta D_{jt} \le D_{jt} \quad \forall j, \forall t \quad (1i)$$

$$0 \le \Delta w_{mt} \le w_{mt} \quad \forall m, \forall t \quad (1j)$$

$$w_{mt} = (w_{mt}^{u} - \hat{w}_{mt})v_{mt}^{u} + (w_{mt}^{l} - \hat{w}_{mt})v_{mt}^{l} + \hat{w}_{mt} \quad (1k)$$

$$\sum_{t=1}^{T}\left(v_{mt}^{u} + v_{mt}^{l}\right) \le \Gamma^{T} \quad \forall m \quad (1l)$$

$$\sum_{m=1}^{M}\left(v_{mt}^{u} + v_{mt}^{l}\right) \le \Gamma^{S} \quad \forall t \quad (1m)$$

$$v_{mt}^{u} + v_{mt}^{l} \le 1 \quad \forall m, \forall t \quad (1n)$$

$$v_{mt}^{u}, v_{mt}^{l} \in \{0,1\} \ \forall m, \forall t \quad (1o)$$

In this formulation, suppose the UC strategy has been given, i.e., $u_{gt}$ is known. Constraint (1b) describes the generation capacity limits of generators. (1c) and (1d) are the ramping rate limits of generators respectively. (1e) is the network power flow limits on transmission lines. (1f) describes the upper and lower limits of the phase angles of nodes and (1g) represents the reference phase angle. (1h) depicts the relaxed power balance requirements for each node with emergency actions including load shedding and wind generation curtailment. (1i) and (1j) give the upper and lower limits of load shedding and wind generation curtailment, respectively. (1k)-(1o) use a polyhedral set to describe the wind generation uncertainty denoted as $W$ as in [14], [15], [17]. Specifically, (1k) depicts the wind generation output; (1l) and (1m) describe the uncertainty budgets over both time horizon and spatial domains, respectively. To determine the parameters in (1k)-(1m), one can refer to [8]. Also, in this model, spinning reserves are provided by the committed units that have been implicitly determined by the given UC strategy. However, additional reserves such as fast-start generators or pumped storage hydro can be incorporated in this model by simply inserting related terms and associated constraints into (1).

The objective function (1a) is the weighted sum of $\Delta w_{mt}$ and $\Delta D_{jt}$, in which the weight coefficients $e_{mt}$ and $f_{jt}$ are non-negative. The value of (1a) is denoted as $F$ and obviously non-negative. $F$=0 means all the wind generation scenarios within $W$ will not cause any load shedding or wind generation curtailment, thus, $W$ is fully admissible. Otherwise, if $F$>0, it means that $W$ is not fully admissible, and some load shedding or wind generation curtailment under certain wind generation scenarios have to be used to recover infeasibility. Therefore, the proposed model can be used to check the admissibility of wind generation.

As $e_{mt}$ and $f_{jt}$ in (1a) are the price of wind generation curtailment and load shedding, respectively, $F$ can be taken as the economic loss in operation [18], inherently providing a criterion for checking the admissibility of wind generation.

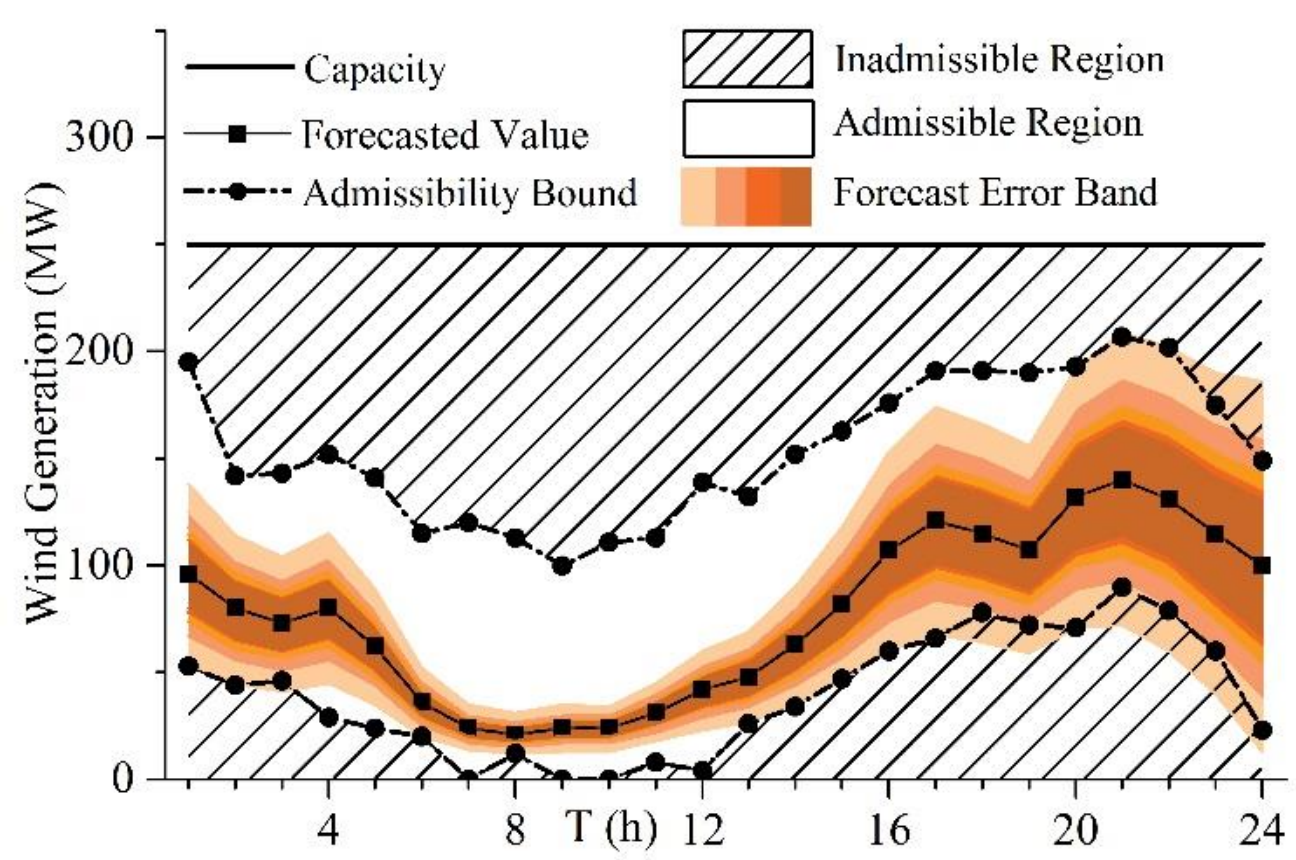

Fig. 1. Schematic diagram of the proposed risk index

### B. A Risk-Based Admissibility Measure

With the proposed model above, ISOs can examine the admissibility of $W$, consequently determine the admissible region of wind generation under the given UC strategy. It is worthy of noting that, besides the UC strategy, the ED strategies also remarkably influence the admissibility of wind generation, since the admissible wind generation can vary among different ED strategies. In this regard, we derive a risk-based measure to make different admissible regions of wind generation comparable in the sense of operational risk brought by inadmissible wind generation. Then, the admissibility assessment problem is converted into an optimization problem, whereby the admissible region of wind generation can be determined reasonably.

As shown in Fig. 1, obviously, the upper and lower boundaries of possible wind generation are the installed capacity of wind farms and zero, respectively. This constitutes the space of wind generation. Assuming the admissibility boundaries of wind generation have been obtained, the space of wind generation can be divided into two parts by the admissibility boundaries (dashed lines with dark dots): the admissible region (unshaded area) and the inadmissible region (shaded area). In the admissible region, no additional emergency regulation is required since any arbitrary realization of wind generation can be fully admitted without breaking the operational feasibility. On the other hand, if the actual wind generation exceeds the admissibility boundaries and enters the inadmissible region, it may lead to undesired power imbalance that cannot be fully handled by the committed units themselves. In such a situation, additional emergency regulations, such as fast-starting units, load shedding or wind generation curtailment may have to be used to recover the operation feasibility. The wind generation forecast error bands under different confidence level $\alpha$ can also be obtained from either the PDF of $\delta_{mt}$ or historical data, as shown in Fig. 1. The part of forecasted wind generation which is within the admissible region is admissible and riskless, while part that is out of the admissible region can lead to operational risk.

As mentioned above, the expected cost for such emergency regulations referred to as operational risk provides an admissibility measure for the wind generation. The operational risk can be calculated by

$$Risk = \sum_{t=1}^{T}\sum_{m=1}^{M}\left( \begin{array}{l} g_t^n \int_{w_{mt}^u-\hat{w}_{mt}}^{w_m^{\max}-\hat{w}_{mt}} \left(\delta_{mt} - w_{mt}^u + \hat{w}_{mt}\right) + \cdots \\ + g_t^p \int_{-\hat{w}_{mt}}^{w_{mt}^l-\hat{w}_{mt}} \left(w_{mt}^l - \delta_{mt} - \hat{w}_{mt}\right) \end{array} \right) y_{mt}(\delta_{mt})\ \mathrm{d}\delta_{mt} \quad (2)$$

where, $w_{mt}^u$ and $w_{mt}^l$ are the upper and the lower boundaries of wind generation, respectively. $g_t^p$ and $g_t^n$ are the prices of the additional positive and negative emergency regulations, respectively. As the proposed risk-based admissibility assessment framework is day-ahead, which means the specific emergency regulations as well as their price are unknown. Therefore, $g_t^p$ and $g_t^n$ can be regarded as the estimation of the real-time price of emergency regulations, which can be derived from operation experiences and historical data. Here we assume the emergency regulations are adequate enough to eliminate any power balance caused by wind power uncertainty. Additionally, the specified emergency regulations are not distinguished in this paper. $\delta_{mt}$ is the wind generation forecast error and $y_{mt}(\cdot)$ is its PDF. In (2), the first and second integral terms represent the operational risk brought by underestimated and overestimated wind generation, respectively. Based on the proposed measure, different admissible regions over multiple periods become comparable.

As mentioned above, the admissible region of wind generation can be determined by minimizing the expected operational risk created by the inadmissible wind generation. This renders a two-stage optimization model as follows.

$$\min_{w_{mt}^u, w_{mt}^l} Risk \quad (3a)$$

$$s.t. \quad \hat{w}_{mt} \le w_{mt}^u \le w_{mt}^{\max} \quad (3b)$$

$$0 \le w_{mt}^l \le \hat{w}_{mt} \quad (3c)$$

$$\begin{array}{l} w_{mt}^u \\ w_{mt}^l \end{array} \in \left\{ \begin{array}{l} \left(\max\limits_{v^u, v^l} \min\limits_{p,\theta,\Delta w,\Delta D} F\right) = 0 \quad (3d) \\ s.t. \quad (1b)\text{-}(1o) \end{array} \right\}$$

where, (3a) represents the operational risk defined by (2); (3b) and (3c) are the constraints for the boundary points, $w_{mt}^u$ and $w_{mt}^l$ , respectively. Mathematically, (3) is a two-stage robust optimization problem, where the first stage problem is (3a) with constraints (3b)-(3c) and the second stage problem is (3d) with constraints (1b)-(1o). The first stage problem ensures the optimality of $w_{mt}^u$ and $w_{mt}^l$ in terms of operational risk, while the second stage problem guarantees that $w_{mt}^u$ and $w_{mt}^l$ satisfy the admissibility criterion.

***Remark:*** if the system operator is a risk seeker, he may allow certain operational loss to some extent to decrease the operational cost for committing additional costly reserve, particularly as the cost for wind generation curtailment is relative low. In this regard, the admissible region will not be riskless anymore and the original admissibility boundary becomes a relaxed one associated with certain allowed loss, which can be regarded as a loss-constrained admissibility boundary (LCAB). LCAB can be computed by solving (3) in which (3d) is replaced by (3e). In (3e), $C_{\text{loss}}$ represents the maximum tolerable operational loss within the LCAB.

$$\left(\max_{v^u, v^l} \min_{p,\theta,\Delta w,\Delta D} F\right) \le C_{loss} \quad (3e)$$

Obviously, (3d) is nothing but a special case of (3e) by just letting $C_{\text{loss}} = 0$.

Note that (3a) is still difficult to calculate due to its nonlinearity and lack of exact distribution information of $\delta_{mt}$. To circumvent this problem, we use the piecewise linear approximation method (PLA) [19] to obtain the approximate linear form of (3a) and rewrite (3) as follows.

$$\min_{w_{mt}^u, w_{mt}^l, Q_{mt}^p, Q_{mt}^n} G = \sum_{t=1}^{T}\sum_{m=1}^{M}(Q_{mt}^p + Q_{mt}^n) \quad (4a)$$

$$s.t. \quad (3b)\text{-}(3c)$$

$$Q_{mt}^p \ge a_{mtsz}^p w_{mt}^u + b_{mtsz}^p \quad \forall m, \forall t, s = 0,1\ldots,S, z = 0,1,\ldots,Z-1. \ (4b)$$

$$Q_{mt}^n \ge a_{mtsz}^n w_{mt}^l + b_{mtsz}^n \quad \forall m, \forall t, s = 0,1\ldots,S, z = 0,1,\ldots,Z-1. \ (4c)$$

$$w_{mt}^u, w_{mt}^l \in \left\{ \begin{array}{c} (3d) \\ s.t. \quad (1b)\text{-}(1o) \end{array} \right\}$$

where, (4a) is a linear approximation of (3a); (4b) and (4c) are auxiliary constraints induced by the PLA treatment. In (4), $a_{mtsz}^p, a_{mtsz}^n, b_{mtsz}^p, b_{mtsz}^n$ are constant coefficients of the piecewise linear approximation; *s* and *z* are ordinal number generated during the PLA treatment; *S* and *Z* are the maximum values of *s* and *z*, respectively. The details of PLA treatment can be found in the Appendix. Obviously, the scales of (4b) and (4c) have remarkable influence on either the computing efficiency or the accuracy of the proposed model, which are chosen according to the desired accuracy and the model scale.

## III. Solution Method

In this section, we will derive the solution method to solve the admissibility assessment problem (4a) with constraints (3b)

-(3d) and (4b)-(4c). Firstly, the admissibility checking sub-problem (3d) with constraints (1b)-(1o) is considered. Its compact form is as follows.

$$\max_{\mathbf{v}} \min_{\mathbf{p},\boldsymbol{\theta},\Delta\mathbf{w},\Delta\mathbf{D}} F^R = \mathbf{e}^T \Delta\mathbf{w} + \mathbf{f}^T \Delta\mathbf{D} \tag{5a}$$

$$s.t. \quad \mathbf{Jp} + \mathbf{K}\boldsymbol{\theta} + \mathbf{M}\Delta\mathbf{w} + \mathbf{N}\Delta\mathbf{D} + \mathbf{P}(\mathbf{w} \circ \mathbf{v}) + \mathbf{Xv} \le \mathbf{r} \tag{5b}$$

$$\mathbf{Sv} \le \mathbf{u} \tag{5c}$$

In (5), $\mathbf{v}$ is the binary vector variable depicting wind generation uncertainty. $\mathbf{p}$ represents the output vector of generators. $\boldsymbol{\theta}$ represents the phase angle vector of nodes. $\Delta\mathbf{w}$ is the wind generation curtailment vector and $\Delta\mathbf{D}$ is the load shedding vector. $\mathbf{J}, \mathbf{K}, \mathbf{M}, \mathbf{N}, \mathbf{P}, \mathbf{X}, \mathbf{S}, \mathbf{e}, \mathbf{f}, \mathbf{r}, \mathbf{u}$ are constant coefficient matrix and can be derived from (1a)-(1o). $\mathbf{w}$ represent the boundary of wind generation output and is fixed in (5). Specially, $\mathbf{w} \circ \mathbf{v}$ is a Hadamard product. (5) is a bi-level mixed integer linear program (MILP) and can be solved by many effective methods, such as the Karush-Kuhn-Tucker (KKT) condition based method [20], and the strong duality theory based method [20], [21]. In this paper, the inner problem is replaced by its dual problem to reformulate (5a) as a single-level bilinear program. It can be solved by either the big-M linearization method [22] or the outer approximation (OA) method [8]. As the big-M linearization method is proved effective with high efficiency and accuracy in practice, this paper adopts it to solve the admissibility checking problem (5). The compact formulation of dual problem of (5) is as follows.

$$\max_{\mathbf{v},\boldsymbol{\lambda}} F^R = \boldsymbol{\lambda}^T \mathbf{r} - \boldsymbol{\lambda}^T \mathbf{Xv} - \boldsymbol{\lambda}^T \mathbf{P}(\mathbf{w} \circ \mathbf{v}) \tag{6a}$$

$$s.t. \quad [\mathbf{J} \vdots \mathbf{K} \vdots \mathbf{M} \vdots \mathbf{N}]^T \boldsymbol{\lambda} \le [\mathbf{0}^T \vdots \mathbf{0}^T \vdots \mathbf{e}^T \vdots \mathbf{f}^T]^T \tag{6b}$$

$$\boldsymbol{\lambda} \le \mathbf{0} \tag{6c}$$

(5c)

where, $\boldsymbol{\lambda}$ is the dual variable vector of inner problem of (5a). Noticed that there is bilinear terms in (6a), auxiliary variables and constraints are introduced to replace them which transfer (6) into a MILP problem as follows.

$$\max_{\mathbf{v},\boldsymbol{\gamma},\boldsymbol{\lambda}} F^R = \boldsymbol{\lambda}^T \mathbf{r} - \boldsymbol{\gamma}^T \mathbf{d} \tag{7a}$$

$s.t.$ (5c), (6b)-(6c)

$$-M_{big} \mathbf{v} \le \boldsymbol{\gamma} \le \mathbf{0} \tag{7b}$$

$$-M_{big}(1 - \mathbf{v}) \le \boldsymbol{\lambda} - \boldsymbol{\gamma} \le \mathbf{0} \tag{7c}$$

where, $\boldsymbol{\gamma}$ is the auxiliary variable vector, $\mathbf{d}$ is a constant vector and can be derived from the following formula.

$$\boldsymbol{\lambda}^T \mathbf{Xv} + \boldsymbol{\lambda}^T \mathbf{P}(\mathbf{w} \circ \mathbf{v}) = \sum_i \sum_y d_{iy} \lambda_i v_y = \boldsymbol{\gamma}^T \mathbf{d}, \quad \gamma_{iy} = \lambda_i v_y \tag{8}$$

(7b) and (7c) are auxiliary constraints generated during objective function linearization using the big-M method. $M_{big}$ is sufficient large positive real number. Thus, (7) result in a standard single-level MILP, which can be solved easily by using commercial solvers such as CPLEX.

So far, the admissibility assessment problem has been reformulated to be (4a) with constraints (3b)-(3c), (4b)-(4c), (5c), (6b)-(6c) and (7a)-(7c), which is a two-stage robust optimization problem. To solve this problem, a key step is to generate proper feasibility cuts based on the solution of the feasibility checking sub-problem (7), and then augment them into the master problem in each iteration. In light of [8], the feasible cut in iteration $k$ can be constructed as

$$-\boldsymbol{\lambda}_k^T \left( \mathbf{P}(\mathbf{w} \circ \mathbf{v}_k) - \mathbf{P}(\mathbf{w}_k \circ \mathbf{v}_k) \right) \le -F_k^R + C_{loss} \tag{9}$$

In (9), $F_k^R$ is the value of (7a) in iteration $k$; $\boldsymbol{\lambda}_k$ and $\mathbf{v}_k$ are the optimal solution of (7) in iteration $k$; $\mathbf{w}_k$ is the optimal solution of (4) in iteration $k$. Note that in the iteration, only feasibility cuts are involved. Motivated by [17], it is desired to use optimality cuts to speed up the convergence. In this regard, a penalty term $K \cdot \eta$ is added into the original objective function (3a). $\eta$ is an auxiliary variable that represents the upper bound of the solution of the feasibility checking sub-problem (7), while $K$ is a large enough positive number. Similar treatment has been used in Bender's decomposition approaches [23]. As a result, the optimality cuts can be generated and the transformed model can be solved efficiently by using the standard C&CG algorithm. For simplicity of the proposed algorithm, the compact form of (4a) with constraints (3b)-(3c) and (4b)-(4c) is written as follows.

$$\min_{\mathbf{w},\mathbf{Q}} G = \mathbf{1}^T \mathbf{Q} \tag{10a}$$

$$s.t. \quad \mathbf{Aw} + \mathbf{BQ} \le \mathbf{c} \tag{10b}$$

where, $\mathbf{Q}$ is the operational risk vector; $\mathbf{A}, \mathbf{B}, \mathbf{c}$ are constant coefficient matrix and can be derived from (3b)-(3c) and (4b)-(4c).

***Algorithm I: Admissibility assessment***

***Step 1:*** *set* $l=0$, $O = \emptyset$, $G_0 = +\infty$.

***Step 2: Solve the following problem and obtain*** $\mathbf{w}_{k+1}$

$$\min_{\mathbf{w},\mathbf{Q},\mathbf{p},\boldsymbol{\theta},\Delta\mathbf{w},\Delta\mathbf{D},\eta} G + K \cdot \eta \tag{11a}$$

***s.t.*** ***(3b)-(3c), (4b)-(4c)***

$$\mathbf{e}^T \Delta\mathbf{w}^k + \mathbf{f}^T \Delta\mathbf{D}^k \le \max\{\eta, C_{loss}\} \quad \forall k \in \mathbf{O} \tag{11b}$$

$$\mathbf{Jp}^k + \mathbf{K}\boldsymbol{\theta}^k + \mathbf{M}\Delta\mathbf{w}^k + \mathbf{N}\Delta\mathbf{D}^k + \mathbf{P}(\mathbf{w} \circ \mathbf{v}^k) + \mathbf{Xv}^k \le \mathbf{r}, \ \forall k \le l \tag{11c}$$

$$-\boldsymbol{\lambda}_k^T \left( \mathbf{P}(\mathbf{w} \circ \mathbf{v}_k) - \mathbf{P}(\mathbf{w}_k \circ \mathbf{v}_k) \right) \le -F_k^R + C_{loss}, \ \forall k \le l \tag{11d}$$

***Step 3: If*** $|G_k - G_{k-1}| < \epsilon$, ***terminate. Otherwise, solve (7), update*** $\eta$, ***derive the optimal solution*** $\mathbf{v}_{k+1}, \boldsymbol{\lambda}_{k+1}$ ***, create variables*** $\mathbf{p}^{k+1}, \boldsymbol{\theta}^{k+1}, \Delta\mathbf{w}^{k+1}, \Delta\mathbf{D}^{k+1}$ ***and add the following constraints***

$$\mathbf{e}^T \Delta\mathbf{w}^{k+1} + \mathbf{f}^T \Delta\mathbf{D}^{k+1} \le \max\{\eta, C_{loss}\} \tag{11e}$$

$$\mathbf{Jp}^{k+1} + \mathbf{K}\boldsymbol{\theta}^{k+1} + \mathbf{M}\Delta\mathbf{w}^{k+1} + \mathbf{N}\Delta\mathbf{D}^{k+1} + \mathbf{P}(\mathbf{w} \circ \mathbf{v}^{k+1}) + \mathbf{Xv}^{k+1} \le \mathbf{r} \tag{11f}$$

$$-\boldsymbol{\lambda}_{k+1}^T \left( \mathbf{P}(\mathbf{w} \circ \mathbf{v}_{k+1}) - \mathbf{P}(\mathbf{w}_{k+1} \circ \mathbf{v}_{k+1}) \right) \le -F_{k+1}^R + C_{loss} \tag{11g}$$

***Update*** $l=l+1$, $\mathbf{O} = \mathbf{O} \cup \{l+1\}$ ***and go to Step 2.***

In Algorithm 1, $\epsilon$ represents the convergence gap. When $\mathbf{w}$ satisfy the admissibility criterion (3d), we must have $\eta = C_{loss}$, which means that (11a) is equivalent to (3a) plus $C_{loss}$ and the correctness of the optimal solution can be ensured. Theoretically, this is guaranteed by imposing constraints (11d) and (11g). Compared with Algorithm 1, the standard C&CG algorithm does not have to consider the feasibility constraints (11d) and (11g) since (5) has a relaxed formulation, which always ensure the feasibility.

## IV. Illustrative Example

In this section, we present numerical experiments carried on the modified IEEE 9-bus system with one wind farm to show the effectiveness of the proposed model and algorithm. The experiments are performed on a PC with Intel(R) Core(TM) 2 Duo 2.2 GHz CPU and 4 GB memory. All algorithms are implemented on MATLAB and programmed using YALMIP. The MILP solver is CPLEX 12.6.

### A. The Modified IEEE 9-bus System

The tested system has 3 generators and 9 transmission lines. A wind farm is connected to the system at bus 1 with an installed capacity of 250 MW. The grid's parameters can be found in [24]. The generators' parameters can be found in Table I. The load curve and the day-ahead forecast of wind generation are both scaled down from the day-ahead curve of California ISO as shown in Fig. 2. The maximum, the minimum, and the average proportion of the wind generation during the day are 29.3%, 5.99%, and 18.7%, respectively. Prices for wind generation curtailment and load shedding can be found in Table II. Assume additional emergency regulation resources can be purchased from the day-ahead reserve market, and the price is listed in Table III. We choose the confidence level $\beta_t = 95\%$, yielding $\Gamma^T \approx 8$ [8]. As only one wind farm is considered in this case, (1m) is unnecessary. In this case, the root mean square error of $\delta_{mt}$ is subject to (12) with $\sigma = 10\%$ and its mean value is zero. In (12), $\sigma$ is a constant parameter.

$$\sigma_{mt} = \sigma \cdot \hat{w}_{mt} \cdot (1 + e^{-(T-t)}) \quad \forall m, \forall t \tag{12}$$

In this case, wind generation forecast error bands is simply derived by Gaussian distribution and as shown in Fig. 2. There are other advanced methods to determine wind generation forecast error bands in the literature, however, it does not influence the computation of admissibility region and beyond the scope of this paper. We choose $\alpha_1^n = 0.5\%$ , $\alpha_2^n = 2.5\%$ , $\alpha_3^n = 49.5\%$ in (13c), which means an eight-piecewise linear distribution approximation is adopted to represent the PDF of $\delta_{mt}$ . We further set $Z$=4 in (4b) and (4c) based on the setting presented in [25]. If not specified, $C_{loss}$ equals zero.

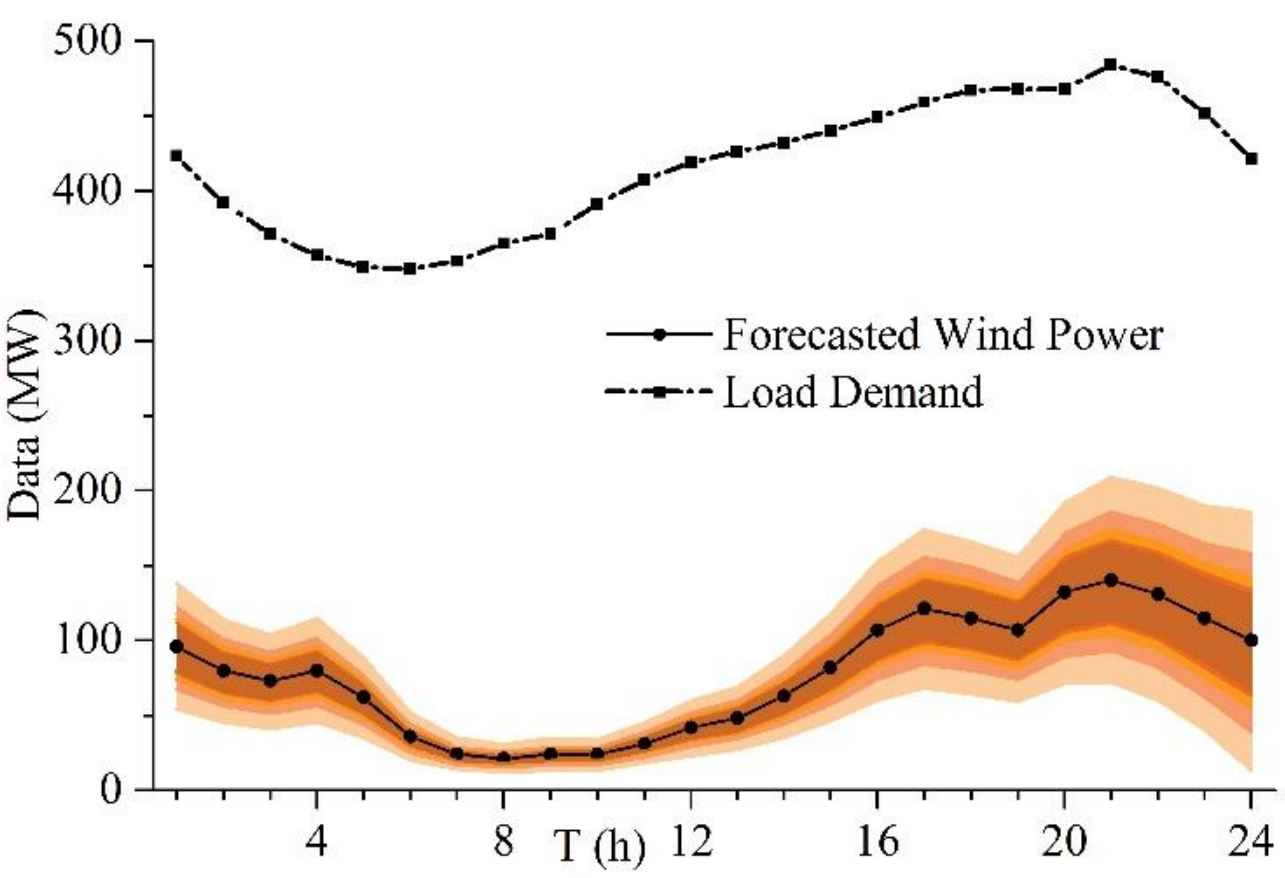

Fig. 2. Data of load demand and forecasted wind power

### B. Effectiveness

In the inter-day operation of a bulk power system, the SRR is usually chosen as $r$=0.03~0.05, depending on the system capacity. In this case, the optimal UC strategy is obtained with $r$=0.05 as the baseline, and then the admissibility assessment of wind generation is conducted using Algorithm I.

The boundaries of the admissible wind generation are shown in Fig. 3 (violet dashed lines). It is found that the confidence levels of the admissible wind generation changes between different periods. During periods 5~14, the confidence levels of admissible wind generation are larger than 99.9%, while in period 17, the confidence level drops to 99%. That is to say, the UC solution with $r$=0.05 cannot guarantee the accommodation of the forecasted wind generation with the given confidence level, i.e., 99.9%. The lowest confidence level is about 95%, which occurs in period 24. Correspondingly, the expected operational risk raised by the inadmissible wind generation is \$41.39. For comparison, the "Do-Not-Exceed" (DNE) limit is computed according to ref. [16], and the result is given in Fig. 3 (green dashed lines). It is shown that, during periods 9~13, the confidence levels of the DNE limit are even smaller than 50%, although the band of the DNE limit is quite large. This is mainly attributed to the fact that the effects of confidence level of wind generation forecast is not considered in the DNE limit. Furthermore, the operational risk beyond the DNE limit is \$3,012, which is much higher than that of admissible region. Operational risk in each period beyond the DNE limit and the admissibility bound is shown in Fig. 4, which is in accordance with the confidence level analysis previously. This indicates that the proposed methodology may provide more reasonable results of the admissibility assessment of wind generation.

TABLE I
PARAMETERS OF GENERATORS

| | $P_{max}$(MW) | $P_{min}$(MW) | $R_+$(MW/h) | $R_-$(MW/h) |
|---|---|---|---|---|
| G1 | 355 | 150 | 50 | 50 |
| G2 | 130 | 20 | 30 | 30 |
| G3 | 55 | 10 | 5 | 5 |

TABLE II
COST COEFFICIENT OF LS AND WGC

| Period | 1-6 | 7-18 | 19-24 |
|---|---|---|---|
| Load shedding (\$/MWh) | 400 | 600 | 500 |
| Wind generation curtailment (\$/MWh) | 40 | 60 | 50 |

TABLE III
RESERVE PRICE IN DAY-AHEAD RESERVE MARKET

| Period | 1-6 | 7-18 | 19-24 |
|---|---|---|---|
| Ramp up (\$/MWh) | 100 | 150 | 125 |
| Ramp down (\$/MWh) | 20 | 30 | 25 |

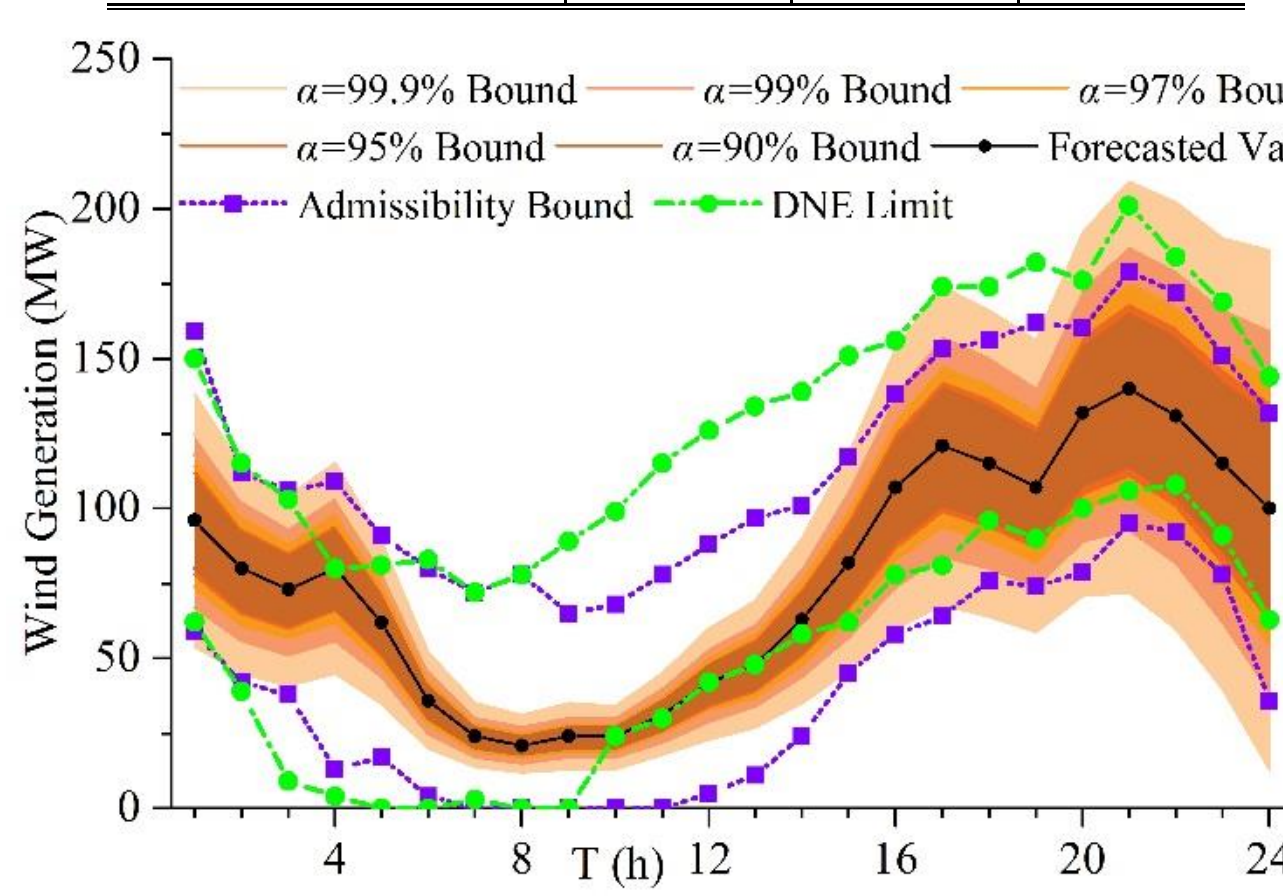

Fig. 3. Confidence levels of wind generation forecast and admissibility boundaries of wind generation.

### C. Computing Efficiency

To investigate the efficiency of the proposed algorithm, we compare the following three algorithms in the test.

**A1**: Algorithm 1.

**A2**: the feasibility cuts based algorithm presented in [8]. It can be implemented by simply setting *K*=0 in Algorithm 1.

**A3**: the standard C&CG algorithm in [17]. It can be implemented by removing (11d) and (11g) from Algorithm 1.

Computational results with different values of *K* are listed in Table IV. **A2** is taken as a baseline for comparison. It can be seen that **A2** can give the correct assessment result of admissibility. However, the computing efficiency seems quite low. When **A3** is employed, it is observed that the same solution can be obtained while the efficiency is much higher when *K* is large enough (*K*=10, *K*=100). Nevertheless, it fails to find the correct solution when *K* is too small (*K*=1). The proposed algorithm **A1** successfully finds the correct solution even more efficiently than the C&CG algorithm no matter *K* is large or small, outperforming its rivals. It also seems that increasing *K* can benefit the efficiency to a large extent. In the following studies, the penalty is identically chosen as *K*=100.

TABLE IV
COMPUTATIONAL RESULTS UNDER DIFFERENT K

| | **A2** | **A3** | | | **A1** | | |
|---|---|---|---|---|---|---|---|
| | | *K*=1 | *K*=10 | *K*=100 | *K*=1 | *K*=10 | *K*=100 |
| Risk($) | 43.19 | 14.24 | 43.19 | 43.19 | 43.19 | 43.19 | 43.19 |
| Iteration | 73 | 57 | 21 | 9 | 33 | 10 | 9 |
| Time(s) | 177.3 | 126.4 | 40.63 | 16.19 | 74.27 | 17.21 | 15.16 |

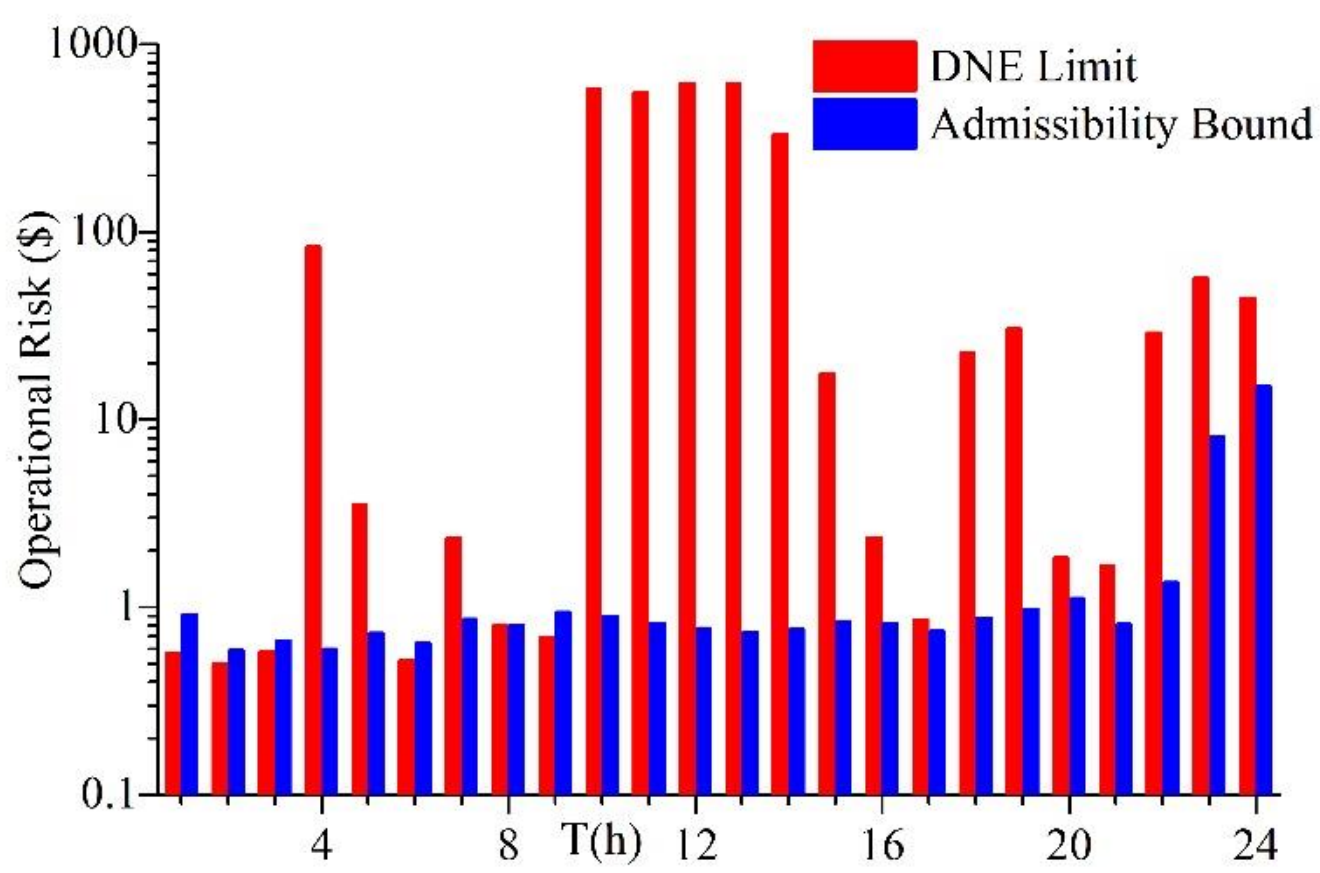

Fig. 4. Operational risk beyond the DNE limit and the admissibility bound.

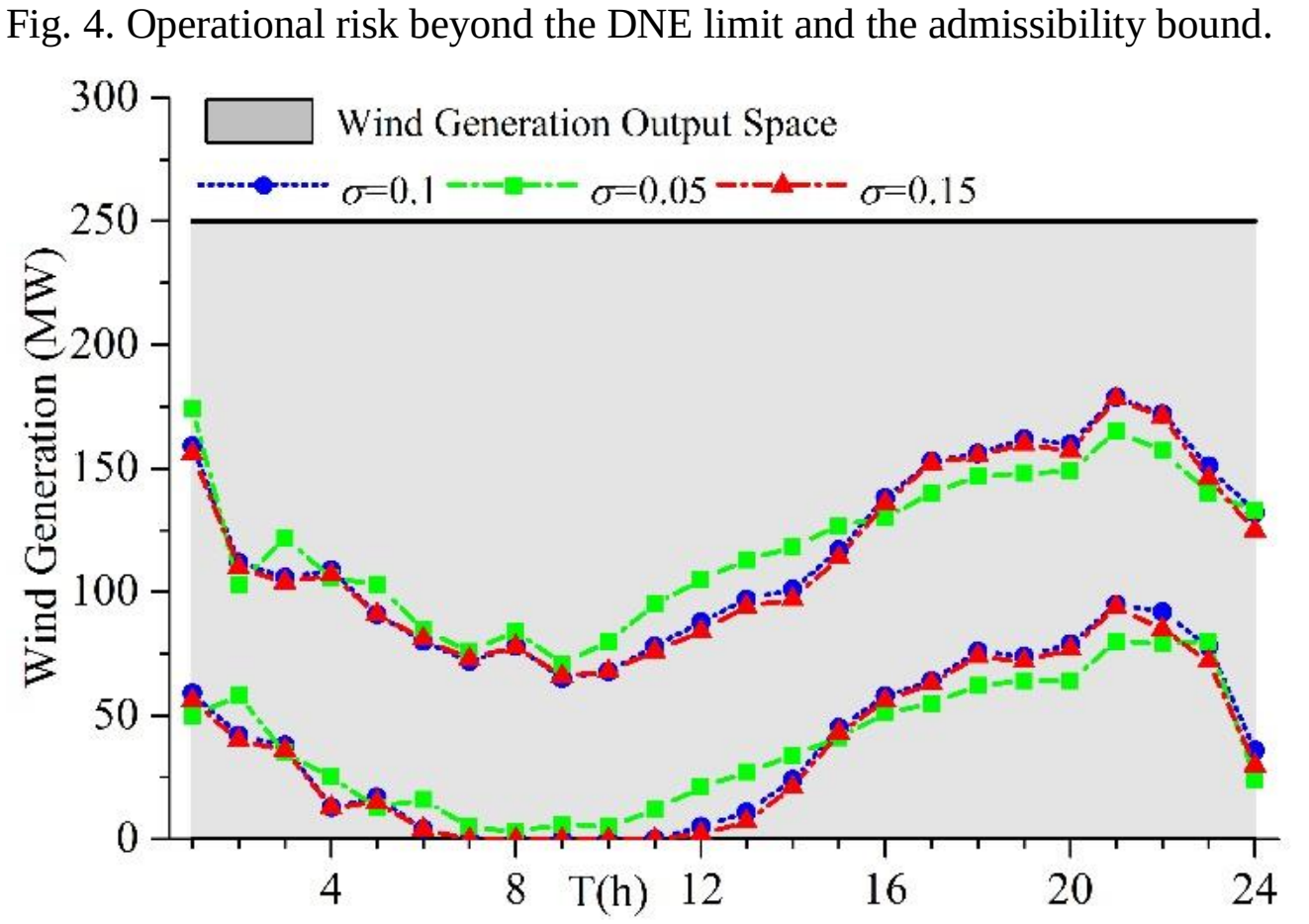

Fig. 5. Admissibility boundaries of wind generation under different σ.

## D. Impacts of Forecast Accuracy

To study the influence of forecast accuracy on the admissibility of wind generation, the admissibility boundaries of wind generation under three levels of forecast errors are computed and shown in Fig. 5. Here, $\sigma = 10\%$ is regarded as the baseline, while $\sigma = 5\%$ and $\sigma = 15\%$ are taken as the high-accuracy and low-accuracy forecast, respectively. It shows that the forecast error can influence the admissibility of wind generation significantly. However, other than shifting the overall boundaries up or down, it changes the admissibility boundaries in an uneven manner. The operational risks in different periods and under different $\sigma$ are shown in Fig. 6. It is found that the operational risks increase along with the growth of forecast error in each of periods. Nevertheless, with the same increment of forecast error, the operational risk changes unevenly in different periods. This result evidently indicates that, the impact of forecast error on the admissibility of wind generation is majorly attributed to the fact that it changes the requirements of ramp reserves rather than the capacity reserve, to cope with the variation and uncertainty of wind generation.

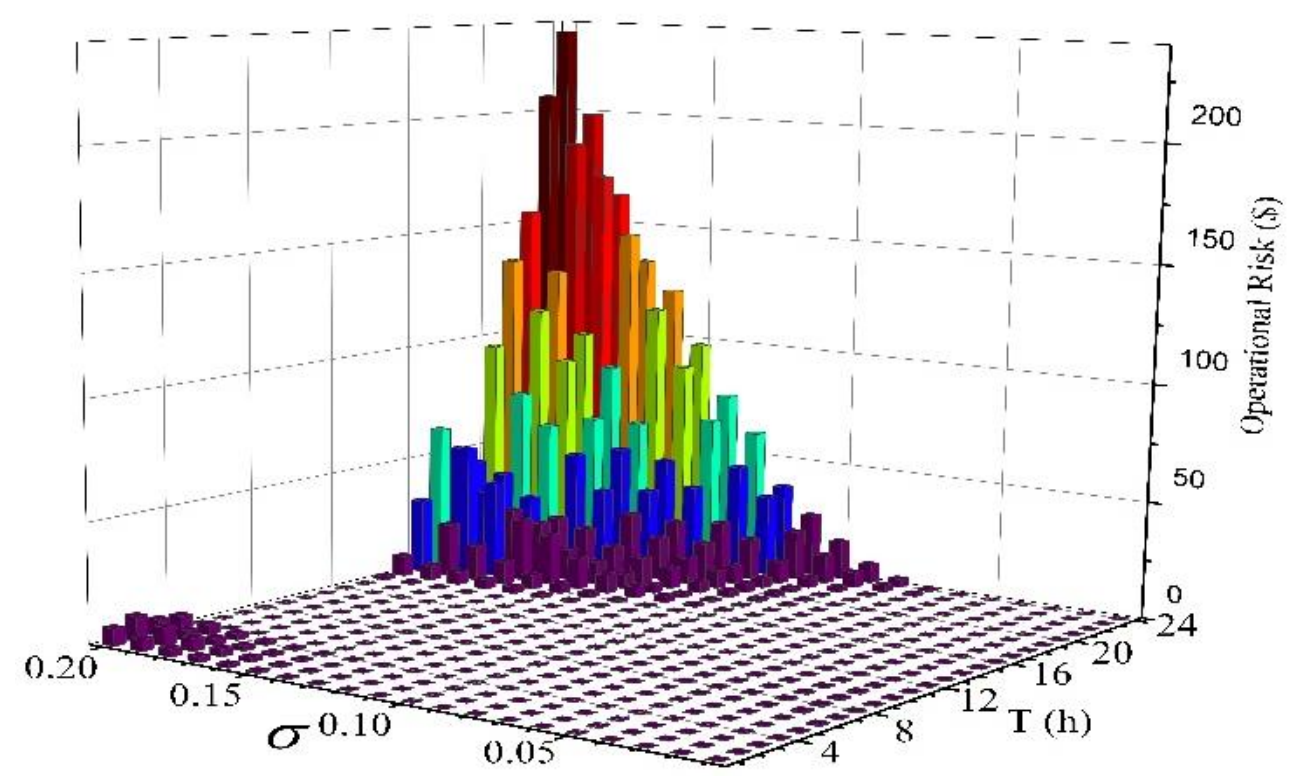

Fig. 6. Risk in different periods and different forecast errors.

As $\sigma$ varies from 0 to 20%, the operational risk is shown in Fig. 7. It is found that the risk almost remains unchanged when $\sigma$ is small enough ($\sigma \leq 6\%$). However, when $\sigma$ becomes large ($\sigma > 6\%$), the operational risk increases approximately exponentially with the growth of $\sigma$. Also, the increment of operational risk when $\sigma$ increases a step (0.01 in this case) is defined as the marginal operational risk, as shown in Fig. 7. This provide the system operator with a quantitative means to evaluate the gain and the cost of improving wind generation forecast accuracy, which may facilitate the management of wind generation forecast.

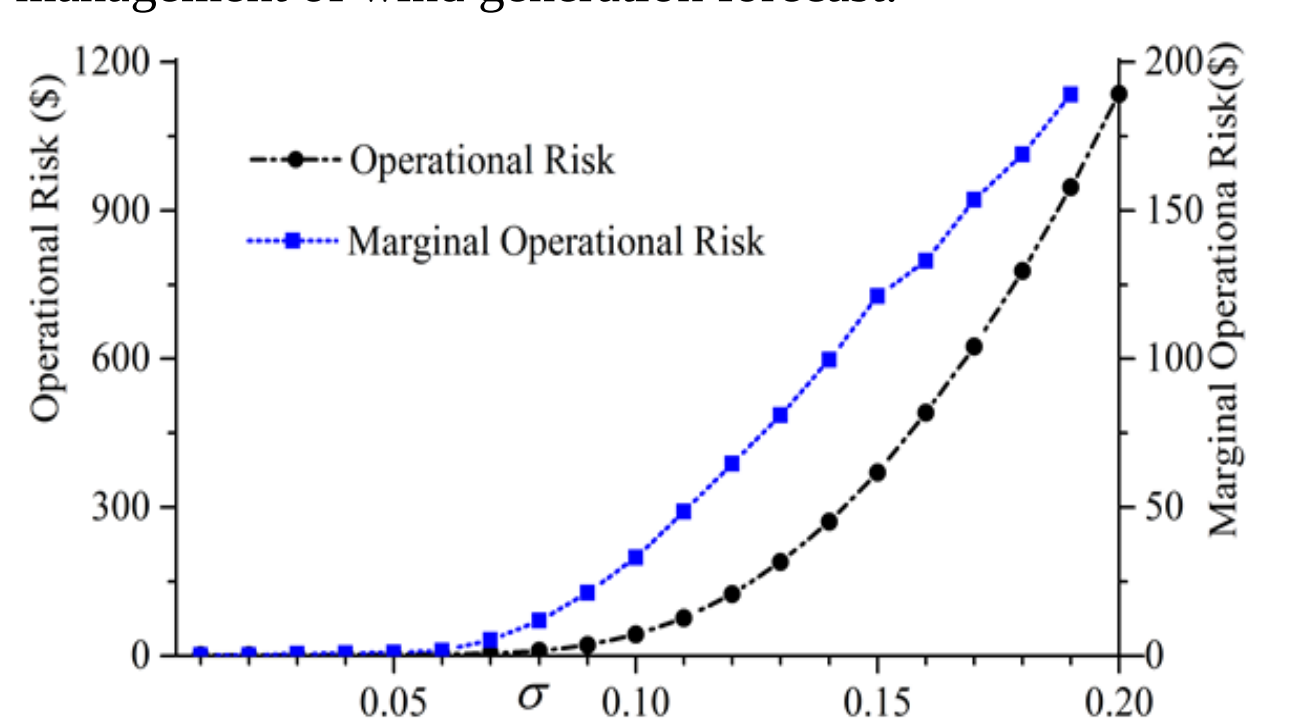

Fig. 7. Operational risk under different σ.

## E. Impacts of $C_{loss}$

To study the influence of $C_{loss}$, admissibility boundaries under three different levels of $C_{loss}$ is computated and shown in Fig. 8. Here, $C_{loss}=0$ is regarded as the base case, while cases with $C_{loss}=1000\$$ and $C_{loss}=2000\$$ are added for comparision. Obviously, the case with $C_{loss}=2000\$$ has the largest admissible region and the case with $C_{loss}=1000\$$ has the second largest admissible region. In period 15-22, the gap among different lower admissibile boundaries are remarkable, due to the price coefficients in Table. III.

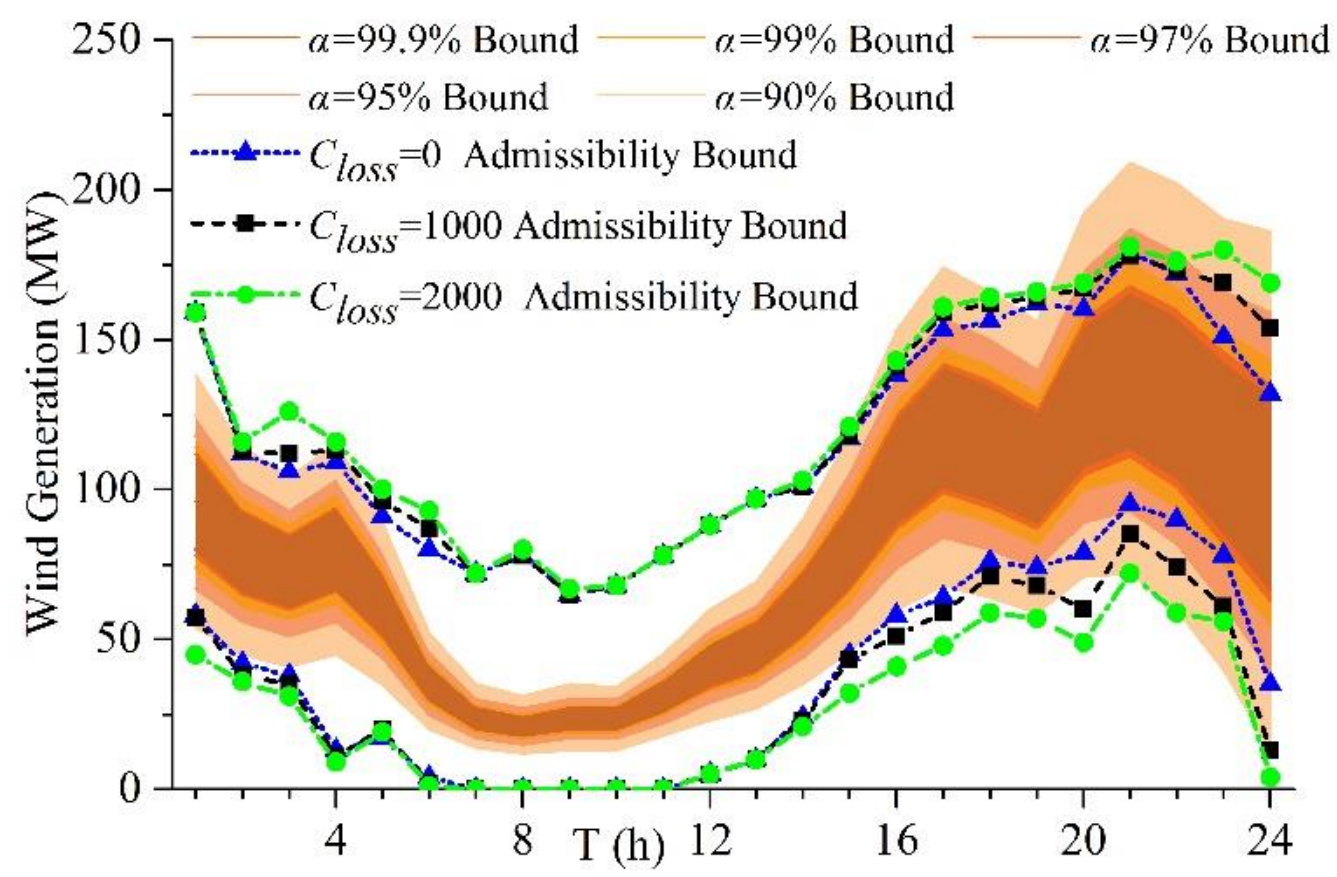

Fig. 8. Admissibility boundaries of wind generation under different $C_{loss}$.

### F. Impacts of Uncertainty Budget

To study the influence of uncertainty budget, admissibility boundaries as well as operational risk under different $\Gamma^T$ are demonstrated in Fig. 9 and Table V, respectively. From Fig. 9 and Table V, the larger $\Gamma^T$ is, the larger operational risk will be, which leads to smaller admissibility region. In this case, the admissibility region as well as operational risk will not change when $\Gamma^T \geq 2$. By choosing appropriate $\Gamma^T$, the robustness and conservatism of the admissibility region can be well balanced.

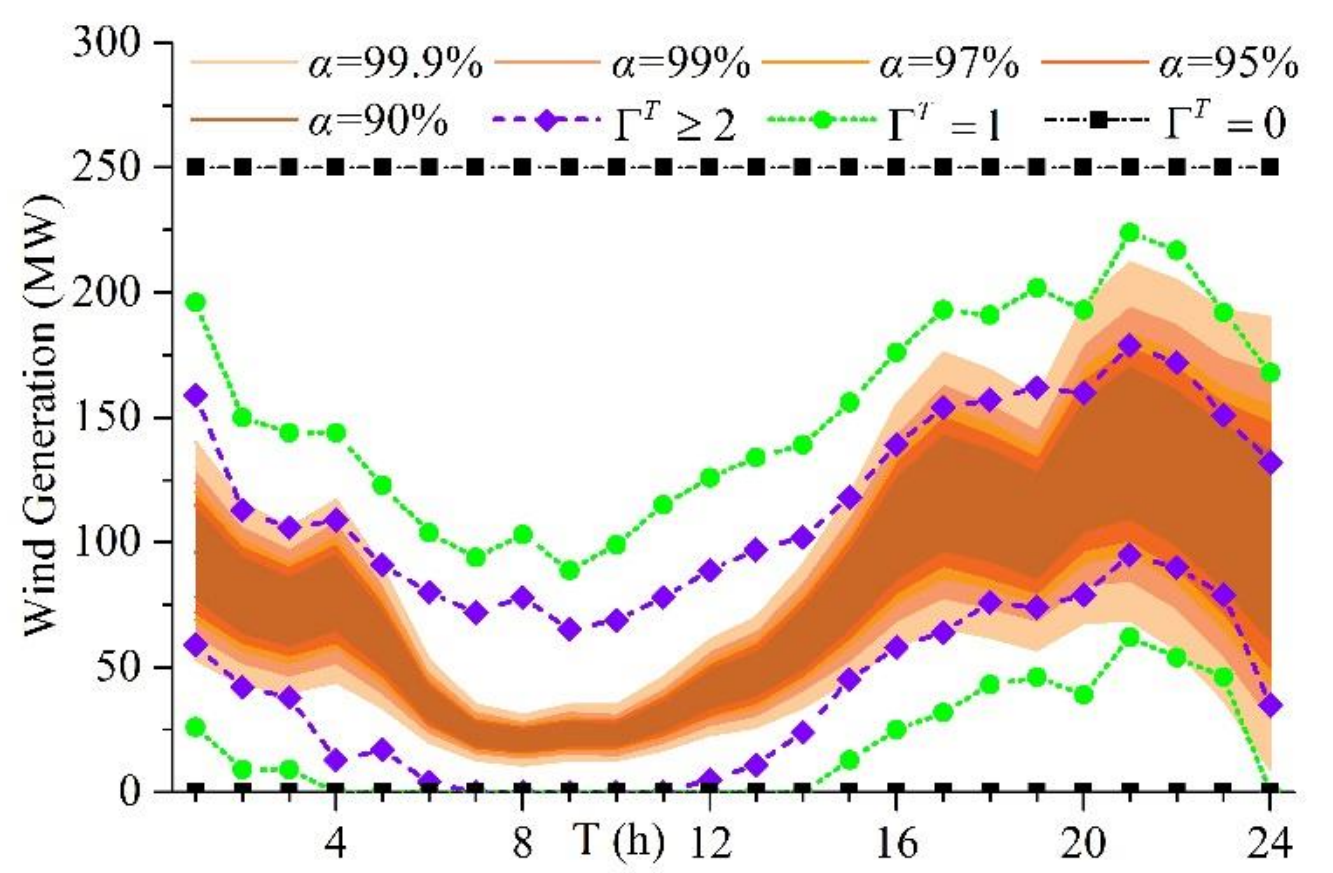

Fig. 9. Admissibility boundaries of wind generation under different $\Gamma^T$.

TABLE V
OPERATIONAL RISK UNDER DIFFERENT $\Gamma^T$

| $\Gamma^T$ | 0 | 1 | $\geq 2$ |
|---|---|---|---|
| Risk ($) | 0 | 7.79 | 43.19 |

## V. A REAL POWER SYSTEM OF CHINA

In this section, we apply the admissibility assessment methodology to the real Guangdong power grid of China. The simulation environment is the same as in section IV.

### A. The Guangdong Power Grid of China

Guangdong power grid has 174 thermal units and 6 wind farms, 453 loads, 1,880 buses and 2,452 transmission lines, among which 100 key transmission lines are selected to be monitored. The total installed capacity of units is 58,744MW. To highlight the influence of wind generation, the installed capacities of six wind farms are modified as follows: 1,500MW in Guangzhou, 1,500MW in Shenzhen, 2,000MW in Dongguan, 2,000MW in Zhanjiang, 2,500MW in Zhuhai, and 2500MW in Foshan, respectively. The total installed capacity of wind generation is about 17% of the installed thermal generation capacity. In particular, we have 24-hour wind generation data of a typical winter day, forecasted nodal loads, and necessary network parameters. The average proportion of integrated wind generation among periods is 21.2%. The price of spinning reserves in the day-ahead reserve market is the same as in section IV. Still, we choose the confidence levels $\beta_t=95\%$ and $\beta_s=95\%$, yielding $\Gamma^T \approx 8$ and $\Gamma^S \approx 4$. The six uncertainty sets of wind generation are formed based on the given confidence levels and parameters. For simplicity, the data of forecasted wind generation and the boundaries of their uncertainty sets are not shown. $\sigma_{mt}$ is the same as in section IV and the mean value of $\delta_{mt}$ is zero. The setting of the PDF approximation of $\delta_{mt}$ is also the same as in section IV. We also set $Z$=4 in (4b) and (4c) according to [25].

### B. Results and Analysis

The admissibility assessment results under different UC strategies are shown in Fig. 10, where the four UC strategies are given by the RUC and the conventional SCUC with reserve rate $r$=0, $r$=0.1 and $r$=0.2 [8], respectively. Specially, the RUC is obtained according to the uncertainty set $W$, where $W$ is constructed from the wind generation prediction result and historical data with the confidence level $\alpha_{mt}=99.9\%$ [6]. It is observed that the admissible region of wind generation under the RUC covers the entire $W$. As the UC decision variables are binary, the admissible wind generation boundaries under RUC are broader than $W$ in certain periods. From the computational results under the three SCUCs with different SRRs, it is found that simply increasing the SRR, is not enough to cope with the variation and uncertainty of wind generation in all periods, although it enlarges the admissible region entirely. Even if $r$ increases up to 0.2, it still cannot guarantee that the admissible region of wind generation covers $W$ with $\alpha_{mt}=90\%$ in certain periods (e.g., periods 1, 2, 22, 23, 24). The main reason is that the conventional SCUC with additional SRR considers the reserve capacity only, while the ramp capacity is not taken into account.

The computational results under the four UCs mentioned above and with different $\sigma$ are shown in Table VI and Fig. 11. It is observed that the operational risk under the RUC is much less sensitive to forecast accuracy compared with the three SCUCs with different SSRs. Fig. 12 shows the operational risks in each period for the four different UCs. It is obvious that the RUC can well hedge against the operational risk brought by the variation and uncertainty of wind generation.

The computing time is listed in Table VI, demonstrating the applicability of our algorithm to realistic power systems.

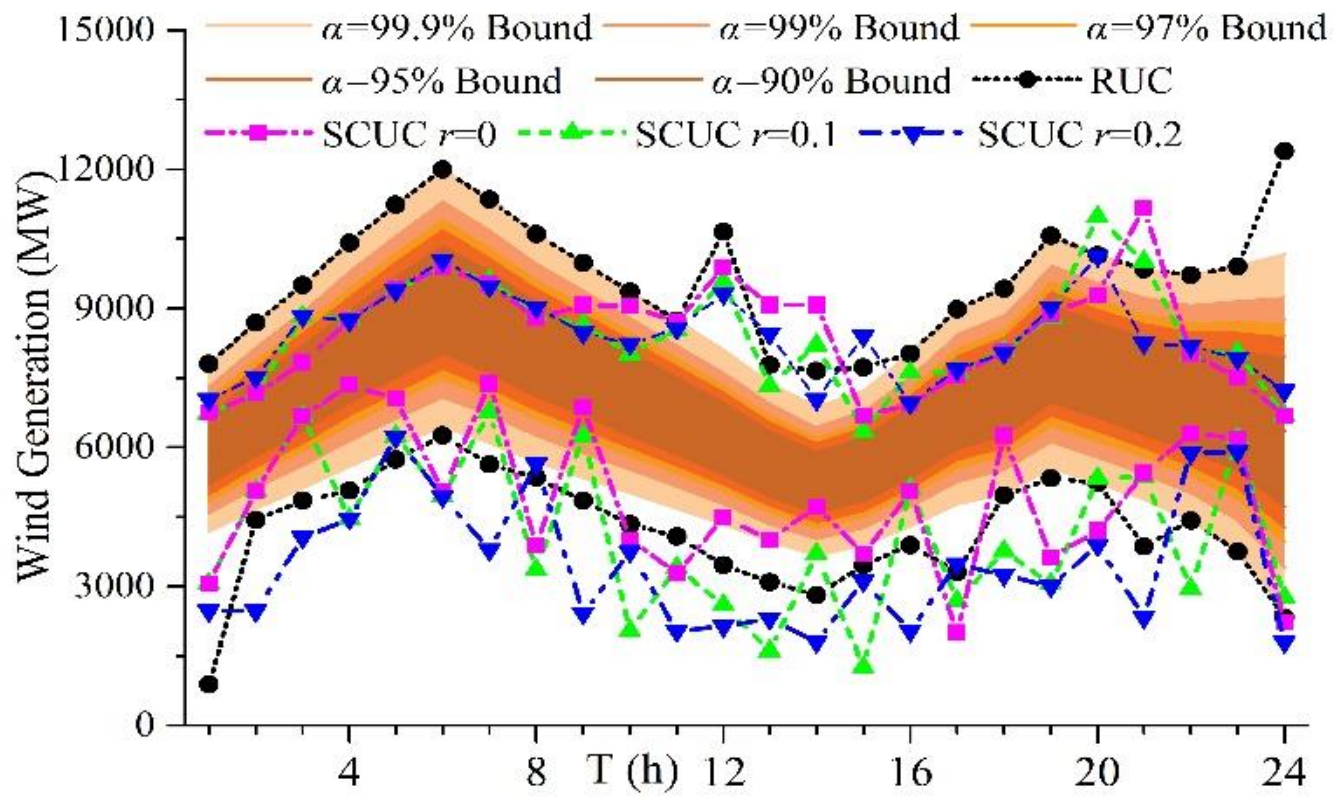

Fig. 10. Admissibility boundaries under different UCs with $\sigma$ =0.1.

TABLE VI
COMPUTATIONAL RESULTS UNDER DIFFERENT UCS AND DIFFERENT σ

| | UC+ED ($) | Risk ($) | | | Time (s) |
|---|---|---|---|---|---|
| | | $\sigma$=0.05 | $\sigma$=0.10 | $\sigma$=0.15 | |
| RUC | $8.32\times10^7$ | $1.04\times10^2$ | $1.11\times10^3$ | $3.93\times10^4$ | 992 |
| SCUC r=0 | $7.32\times10^7$ | $2.12\times10^4$ | $1.65\times10^5$ | $4.34\times10^5$ | 317 |
| SCUC r=0.1 | $7.58\times10^7$ | $1.36\times10^4$ | $9.22\times10^4$ | $2.76\times10^5$ | 305 |
| SCUC r=0.2 | $8.08\times10^7$ | $3.40\times10^3$ | $4.11\times10^4$ | $1.39\times10^5$ | 309 |

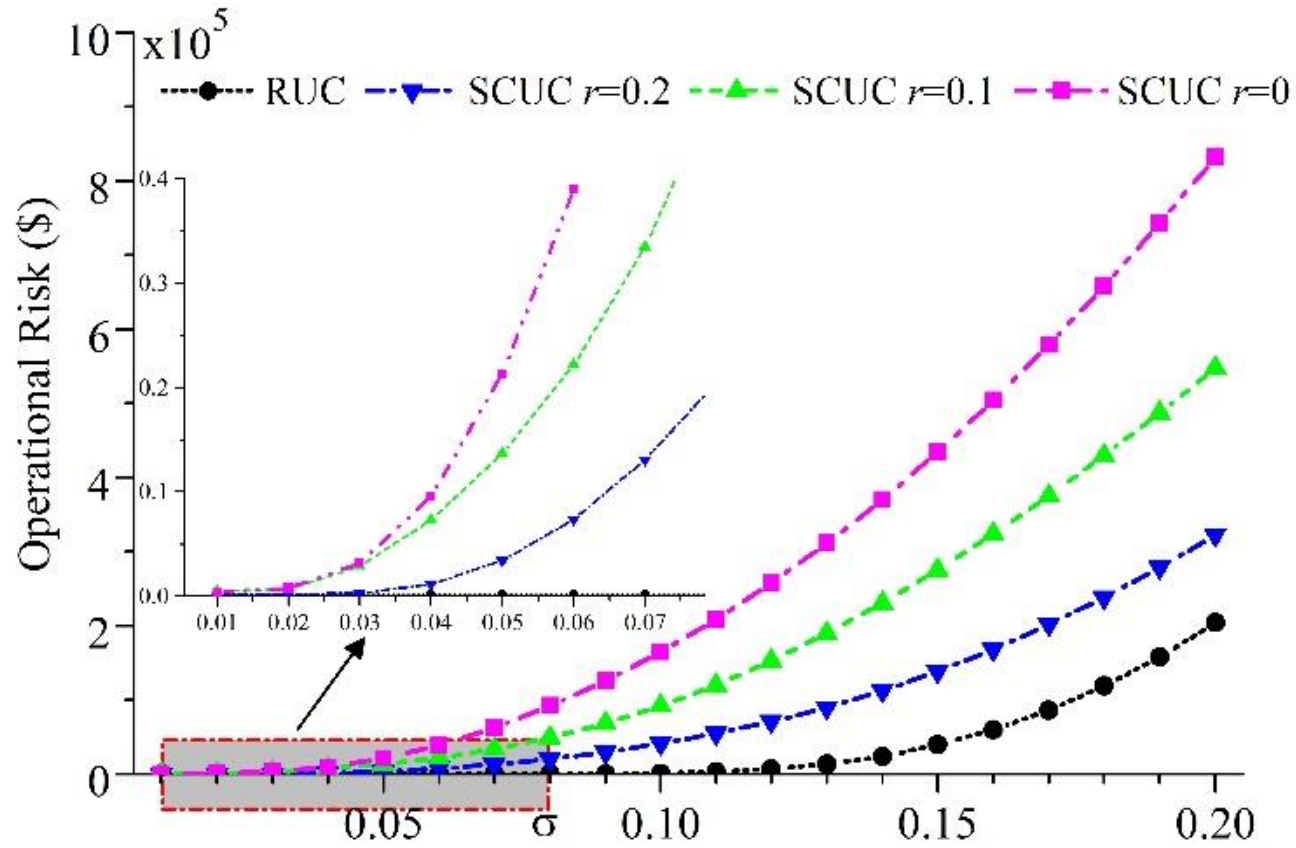

Fig. 11. Operational risks under different UCs and with different σ.

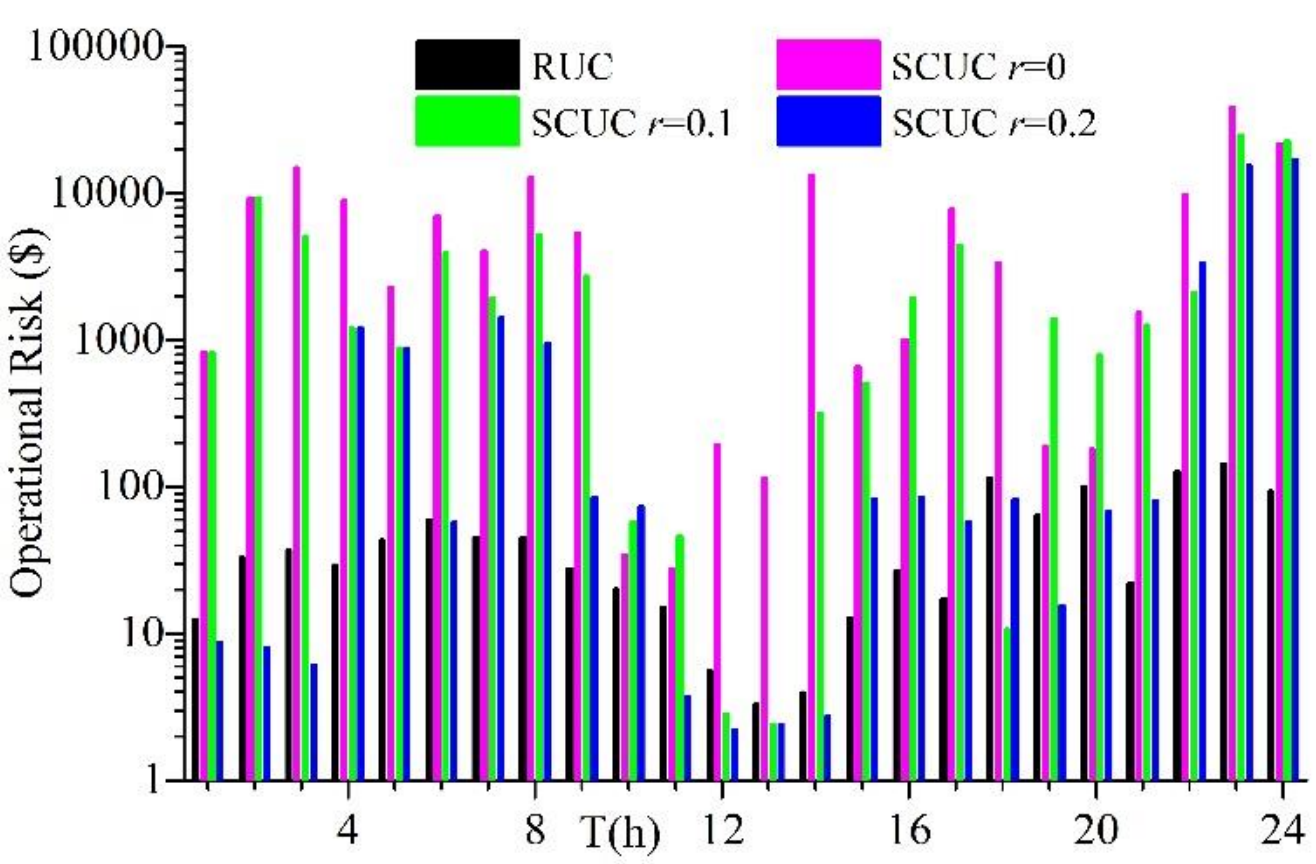

Fig. 12. Operational risks under different UCs and in different periods

## VI. CONCLUSION

In this paper, a risk-based model is proposed for assessing the admissibility of volatile and uncertain wind generation under a given UC strategy. The expected operational risk brought by variation and uncertainty of wind generation is introduced as an admissibility measure of wind generation. As the operational risk depends on the distribution of wind generation forecast error and is nonlinear, a practical linear approximation of the measure is derived. Consequently, a risk-minimization model is derived to characterize the admissible region of wind generation. An iterative algorithm is proposed based on the C&CG algorithm to solve this proposed tri-level model efficiently. Simulations are carried out on the modified IEEE 9-bus system to illustrate the effectiveness of the proposed model and algorithm. It also reveals that the influence of the forecast error of wind generation on the admissibility is mainly attributed to the fact that it increases the requirement of ramp reserve rather than that of capacity reserve. The proposed methodology is also applied to the real Guangdong Power Grid to analyze the admissibility of wind generation as well as the operational risk over the time horizon, demonstrating the practicality of our methodology.

## APPENDIX

In this section, the details of PLA treatment is demonstrated as follows. Firstly, we use the following piecewise linear approximation of $y_{mt}(\cdot)$ to replace its exact form:

$$\sum_{t=1}^{T}\sum_{m=1}^{M}\left(\begin{array}{l} g_t^n \left(\begin{array}{l}\sum_{s^p=0}^{S^p}\int_{Y_{mt}^{-1}(\alpha^p_{s^p+1})}^{Y_{mt}^{-1}(\alpha^p_{s^p})}\left(c^p_{mts^p}\delta_{mt}+d^p_{mts^p}\right)+\\ \int_{w^u_{mt}-\hat{w}_{mt}}^{Y_{mt}^{-1}(\alpha^p_{S^p+1})}\left(c^p_{mt(S^p+1)}\delta_{mt}+d^p_{mt(S^p+1)}\right)\end{array}\right)\left(\delta_{mt}+\hat{w}_{mt}-w^u_{mt}\right) \\ -g_t^p\left(\begin{array}{l}\sum_{s^n=0}^{S^n}\int_{Y^{-1}_{mt}(\alpha^n_{s^n})}^{Y^{-1}_{mt}(\alpha^n_{s^n+1})}\left(c^n_{mts^n}\delta_{mt}+d^n_{mts^n}\right)+\\ \int_{Y^{-1}_{mt}(\alpha^n_{S^n+1})}^{w^l_{mt}-\hat{w}_{mt}}\left(c^n_{mt(S^n+1)}\delta_{mt}+d^n_{mt(S^n+1)}\right)\end{array}\right)\left(\delta_{mt}+\hat{w}_{mt}-w^l_{mt}\right)\end{array}\right)\mathrm{d}\delta_{mt} \quad (13a)$$

where,

$$\begin{array}{l} Y_{mt}^{-1}(\alpha^p_{S^p+2})+\hat{w}_{mt} < w^u_{mt} \le Y_{mt}^{-1}(\alpha^p_{S^p+1})+\hat{w}_{mt},\ \forall m,\forall t \\ Y_{mt}^{-1}(\alpha^n_{S^n+1})+\hat{w}_{mt} < w^l_{mt} \le Y_{mt}^{-1}(\alpha^n_{S^n+2})+\hat{w}_{mt},\ \forall m,\forall t \end{array} \quad (13b)$$

$$\begin{array}{l} \alpha_0^n < \alpha_1^n < \cdots < \alpha_S^n < Y_{mt}(0) \\ \alpha_0^n = 0,\ \ \alpha_s^n + \alpha_s^p = 1,\ \ s = 0,1,...,S. \end{array} \quad (13c)$$

In (13a), $Y_{mt}^{-1}(\cdot)$ represents the inverse function of cumulative density function (CDF) of $\delta_{mt}$ and can be obtained by curve fitting using the historical data. $\alpha_0^n, \alpha_0^p, \alpha_1^n, \alpha_1^p, \dots, \alpha_S^n, \alpha_S^p$ are a set of selected confidence levels of $y_{mt}(\cdot)$ and subject to (13c). In (13c), $Y_{mt}(\cdot)$ is the CDF of $\delta_{mt}$ and can also be obtained by curve fitting using the historical data. $c_{mts}^p, d_{mts}^p, c_{mts}^n, d_{mts}^n$ are subject to (13d). A schematic diagram is shown in Fig. 13.

$$\begin{array}{l} c^n_{mts} = \begin{cases} -\left(y_{mt}(0) - y_{mt}\left(Y^{-1}_{mt}(\alpha^n_S)\right)\right)/Y^{-1}_{mt}(\alpha^n_S), & s = S. \\ \left(y_{mt}\left(Y^{-1}_{mt}(\alpha^n_{s+1})\right) - y_{mt}\left(Y^{-1}_{mt}(\alpha^n_s)\right)\right)/\left(Y^{-1}_{mt}(\alpha^n_{s+1}) - Y^{-1}_{mt}(\alpha^n_s)\right), & other \end{cases} \\ d^n_{mts} = y_{mt}\left(Y^{-1}_{mt}(\alpha^n_s)\right) - c^n_{mts}Y^{-1}_{mt}(\alpha^n_s),\ \forall s \\ c^p_{mts} = \begin{cases} -\left(y_{mt}(0) - y_{mt}\left(Y^{-1}_{mt}(\alpha^p_S)\right)\right)/Y^{-1}_{mt}(\alpha^p_S), & s = S. \\ \left(y_{mt}\left(Y^{-1}_{mt}(\alpha^p_{s+1})\right) - y_{mt}\left(Y^{-1}_{mt}(\alpha^p_s)\right)\right)/\left(Y^{-1}_{mt}(\alpha^p_{s+1}) - Y^{-1}_{mt}(\alpha^p_s)\right), & other \end{cases} \\ d^p_{mts} = y_{mt}\left(Y^{-1}_{mt}(\alpha^p_{s-1})\right) - c^p_{mts}Y^{-1}_{mt}(\alpha^p_{s-1}),\ \forall s \end{array} \quad (13d)$$

Formula (13a) is the sum of integration of the piecewise quadratic function that is equivalent to a continuous cubic function. The cubic function can be further linearized by the PLA method. Finally, (2) can be further approximated by a linear expression with auxiliary variables and constraints using PLA technique again, yielding the following expression:

$$Risk = \min_{Q^p_{mt}, Q^n_{mt}} \sum_{t=1}^{T}\sum_{m=1}^{M}(Q^p_{mt} + Q^n_{mt}) \tag{13e}$$

(4b)-(4c)

The expression of $a^p_{mtsz}, a^n_{mtsz}, b^p_{mtsz}, b^n_{mtsz}$ in (4b) and (4c) are as follows. The relationship between $w^u_{mt}$ and $Q^u_{mt}$ is shown in Fig. 14, as a schematic diagram. The relationship between $w^l_{mt}$ and $Q^l_{mt}$ is similar to that of $w^u_{mt}$ and $Q^u_{mt}$, which will not be repeated here.

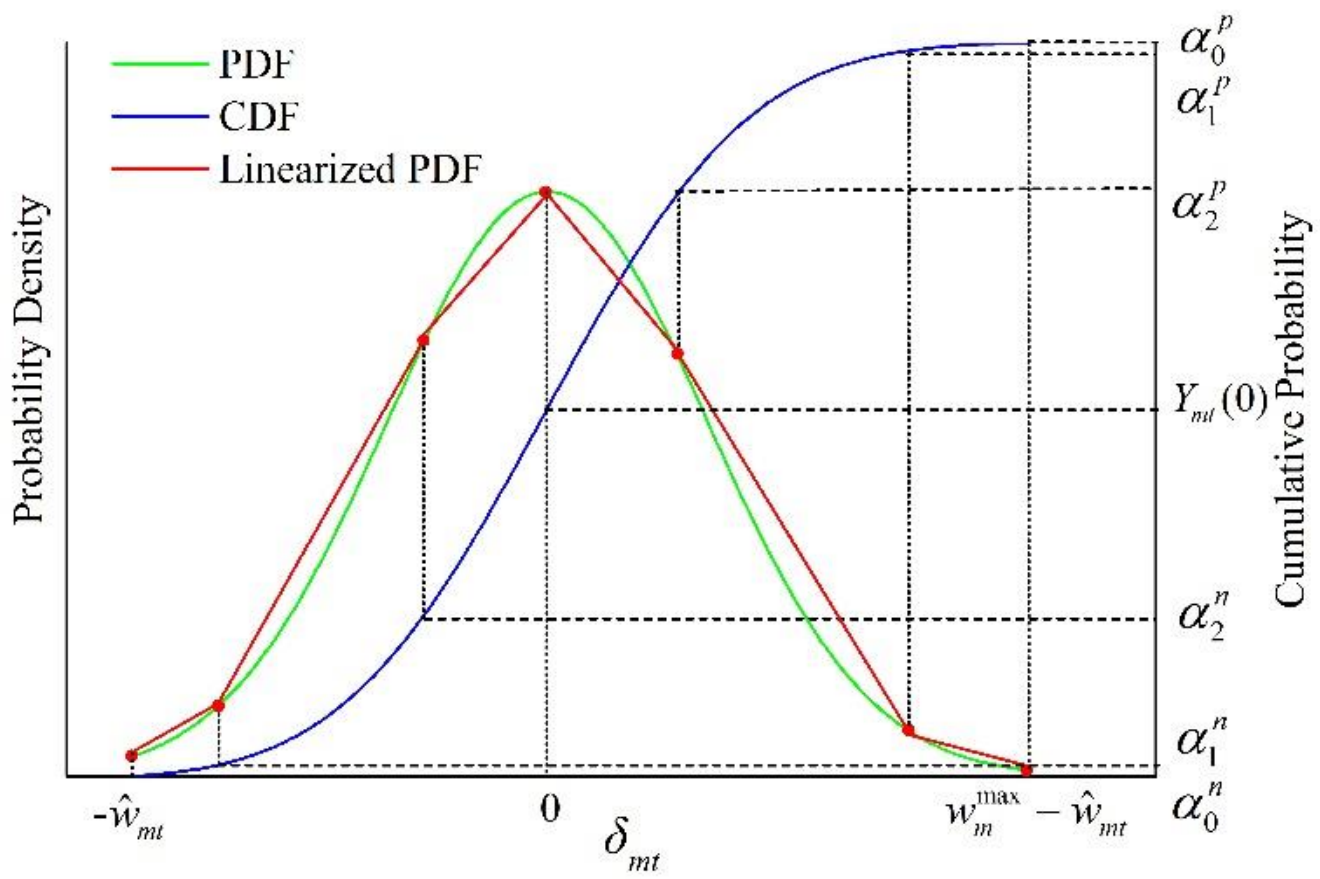

Fig. 13. Piecewise linear PDF of $\delta_{mt}$.

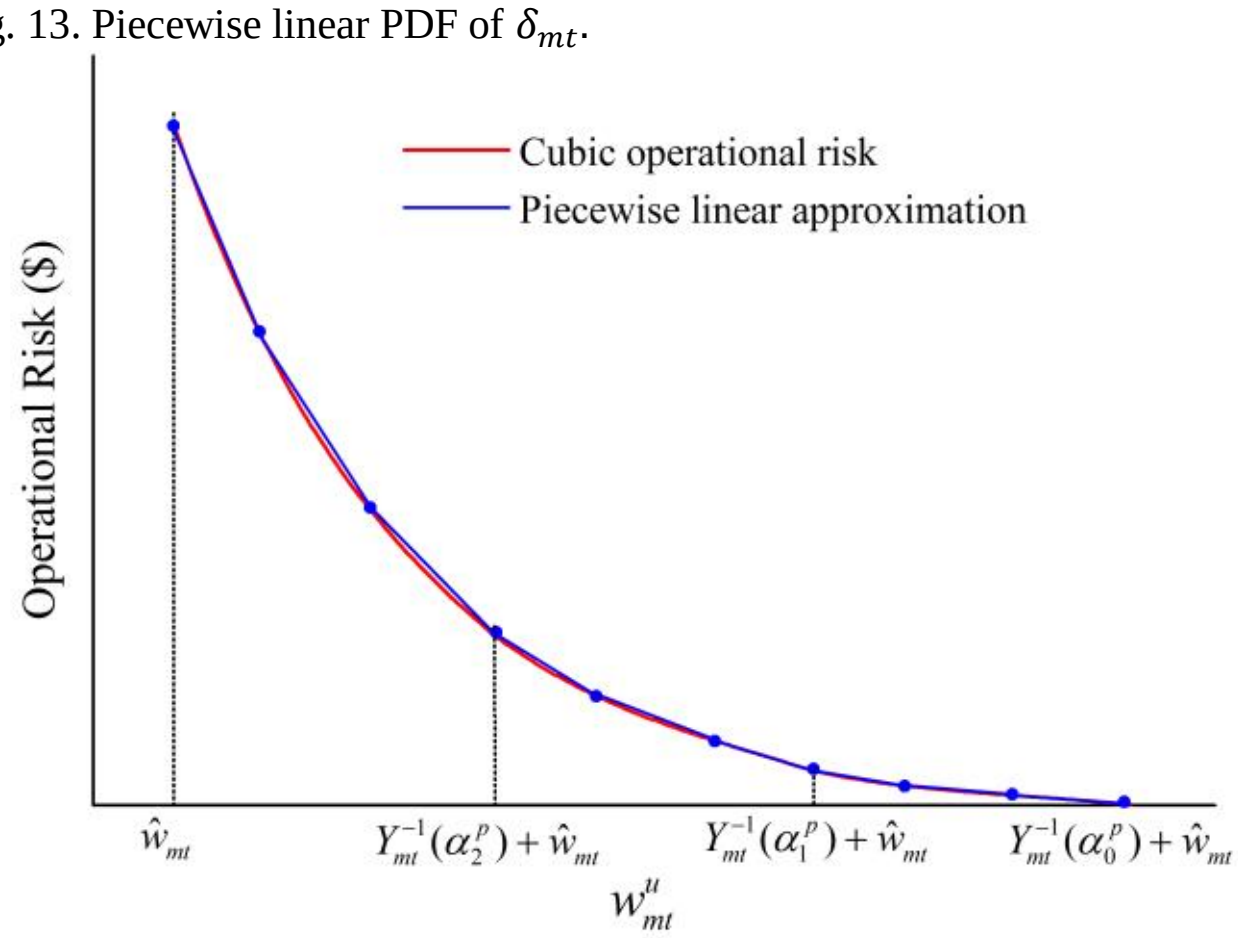

Fig. 14. Piecewise linear relationship between $w^u_{mt}$ and $Q^u_{mt}$.

$$\begin{aligned} a^p_{mtsz} &= g^n_t\left(e^p_{mts(z+1)} - e^p_{mtsz}\right) / \left(f^p_{mts(z+1)} - f^p_{mtsz}\right) \\ b^p_{mtsz} &= g^n_t\left(e^p_{mtsz} - a^p_{mtsz} f^p_{mtsz}\right) \quad \forall m, \forall t, \forall s, z = 0, \cdots, Z-1. \end{aligned} \tag{13f}$$

$$\begin{aligned} a^n_{mtsz} &= g^p_t\left(e^n_{mts(z+1)} - e^n_{mtsz}\right) / \left(f^n_{mts(z+1)} - f^n_{mtsz}\right) \\ b^n_{mtsz} &= g^p_t\left(e^n_{mtsz} - a^n_{mtsz} f^n_{mtsz}\right) \quad \forall m, \forall t, \forall s, z = 0, \cdots, Z-1. \end{aligned} \tag{13g}$$

where,

$$f^p_{mtsz} = Y^{-1}_{mt}(\alpha^p_s) + \hat{w}_{mt} + z\left(Y^{-1}_{mt}(\alpha^p_{s+1}) - Y^{-1}_{mt}(\alpha^p_s)\right) / Z, \forall m, \forall t, \forall s, \forall z.$$

$$e^p_{mtsz} = \begin{cases} c^p_{mts}\left(\left(Y^{-1}_{mt}(\alpha^p_s)\right)^3 - \left(f^p_{mtsz} - \hat{w}_{mt}\right)^3\right)/3 \cdots \\ +\left(c^p_{mts}\hat{w}_{mt} + d^p_{mts} - c^p_{mts} f^p_{mtsz}\right)\left(\left(Y^{-1}_{mt}(\alpha^p_s)\right)^2 - \left(f^p_{mtsz} - \hat{w}_{mt}\right)^2\right)/2 \cdots \\ +\left(d^p_{mts}\hat{w}_{mt} - d^p_{mts} f^p_{mtsz}\right)\left(Y^{-1}_{mt}(\alpha^p_s) - f^p_{mtsz} + \hat{w}_{mt}\right), \forall m, \forall t, \forall z, s = 0. \\ c^p_{mts}\left(\left(Y^{-1}_{mt}(\alpha^p_s)\right)^3 - \left(f^p_{mtsz} - \hat{w}_{mt}\right)^3\right)/3 \cdots \\ +\left(c^p_{mts}\hat{w}_{mt} + d^p_{mts} - c^p_{mts} f^p_{mtsz}\right)\left(\left(Y^{-1}_{mt}(\alpha^p_s)\right)^2 - \left(f^p_{mtsz} - \hat{w}_{mt}\right)^2\right)/2 \cdots \\ +\left(d^p_{mts}\hat{w}_{mt} - d^p_{mts} f^p_{mtsz}\right)\left(Y^{-1}_{mt}(\alpha^p_s) - f^p_{mtsz} + \hat{w}_{mt}\right) + \\ +\sum_{k=1}^{s}\left(c^p_{mtk}\left(\left(Y^{-1}_{mt}(\alpha^p_{k-1})\right)^3 - \left(Y^{-1}_{mt}(\alpha^p_k)\right)^3\right)/3 \cdots\right. \\ +\left(c^p_{mtk}\hat{w}_{mt} + d^p_{mtk} - c^p_{mtk} f^p_{mtsz}\right)\left(\left(Y^{-1}_{mt}(\alpha^p_{k-1})\right)^2 - \left(Y^{-1}_{mt}(\alpha^p_k)\right)^2\right)/2 \cdots \\ \left.+\left(d^p_{mtk}\hat{w}_{mt} - d^p_{mtk} f^p_{mtsz}\right)\left(Y^{-1}_{mt}(\alpha^p_{k-1}) - Y^{-1}_{mt}(\alpha^p_k)\right)\right), \forall m, \forall t, \forall z, s \neq 0. \end{cases} \tag{13h}$$

$$f^n_{mtsz} = Y^{-1}_{mt}(\alpha^n_s) + \hat{w}_{mt} + z\left(Y^{-1}_{mt}(\alpha^n_{s+1}) - Y^{-1}_{mt}(\alpha^n_s)\right) / Z \quad \forall m, \forall t, \forall s, \forall z.$$

$$e^n_{mtsz} = \begin{cases} c^n_{mts}\left(\left(Y^{-1}_{mt}(\alpha^n_s)\right)^3 - \left(f^n_{mtsz} - \hat{w}_{mt}\right)^3\right)/3 \cdots \\ +\left(c^n_{mts}\hat{w}_{mt} + d^n_{mts} - c^n_{mts} f^n_{mtsz}\right)\left(\left(Y^{-1}_{mt}(\alpha^n_s)\right)^2 - \left(f^n_{mtsz} - \hat{w}_{mt}\right)^2\right)/2 \cdots \\ +\left(d^n_{mts}\hat{w}_{mt} - d^n_{mts} f^n_{mtsz}\right)\left(Y^{-1}_{mt}(\alpha^n_s) - f^n_{mtsz} + \hat{w}_{mt}\right), \forall m, \forall t, \forall z, s = 0. \\ c^n_{mts}\left(\left(Y^{-1}_{mt}(\alpha^n_s)\right)^3 - \left(f^n_{mtsz} - \hat{w}_{mt}\right)^3\right)/3 \cdots \\ +\left(c^n_{mts}\hat{w}_{mt} + d^n_{mts} - c^n_{mts} f^n_{mtsz}\right)\left(\left(Y^{-1}_{mt}(\alpha^n_s)\right)^2 - \left(f^n_{mtsz} - \hat{w}_{mt}\right)^2\right)/2 \cdots \\ +\left(d^n_{mts}\hat{w}_{mt} - d^n_{mts} f^n_{mtsz}\right)\left(Y^{-1}_{mt}(\alpha^n_s) - f^n_{mtsz} + \hat{w}_{mt}\right) \cdots \\ +\sum_{k}^{s}\left(c^n_{mtk}\left(\left(Y^{-1}_{mt}(\alpha^n_{k-1})\right)^3 - \left(Y^{-1}_{mt}(\alpha^n_k)\right)^3\right)/3 \cdots\right. \\ +\left(c^n_{mtk}\hat{w}_{mt} + d^n_{mtk} - c^n_{mtk} f^n_{mtsz}\right)\left(\left(Y^{-1}_{mt}(\alpha^n_{k-1})\right)^2 - \left(Y^{-1}_{mt}(\alpha^n_k)\right)^2\right)/2 \cdots \\ \left.+\left(d^n_{mtk}\hat{w}_{mt} - d^n_{mtk} f^n_{mtsz}\right)\left(Y^{-1}_{mt}(\alpha^n_{k-1}) - Y^{-1}_{mt}(\alpha^n_k)\right)\right), \forall m, \forall t, \forall z, s \neq 0. \end{cases} \tag{13i}$$